\documentclass[a4paper,12pt]{article}
\usepackage[frenchb]{babel}
\usepackage{amssymb,xypic,latexsym,graphicx,caption2}
\usepackage[matrix,arrow]{xy}
\newtheorem{theoreme}{Th\'eor\`eme}
\newtheorem{cor}{Corollaire} 

\newtheorem{lemme}{Lemme}
\newtheorem{sous-lemme}{Sous-lemme}
\newtheorem{prop}{Proposition}

\newcommand{\QQ}{{\mathbb Q}}
\newcommand{\NN}{{\mathbb N}}
\newcommand{\cucu}{\hfill $\Box$}

\newcommand{\HH}{{\rm H}}
\newcommand{\hh}{{\rm h}}

\newcommand{\RR}{{\mathbb R}}
\newcommand{\C}{{\mathbb C}}
\newcommand{\PP}{{\mathbb P}}

\newcommand{\im}{{\text{Im}}\,}

\newcommand{\codim}{\text{codim}\,}
\newcommand{\kker}{\text{Ker}\,}
\newcommand{\coker}{\text{Coker}\,}
\newcommand{\anneau}{{\cal O}}
\newcommand{\preuve}{\noindent{\it Preuve. --- \ }}
\newcommand{\prim}{{\rm prim}}
\newcommand{\ev}{{\rm ev}}

\newcommand{\tube}{{\rm tube}}
\newcommand{\bord}{{\rm bord}}
\newcommand{\ouv}{{\rm ouv}}
\newcommand{\rel}{{\rm rel}}

\begin{document}
\def\figurename{{\sc Dessin}}
\title{Monodromie d'une famille d'hypersurfaces}
\author{Ania Otwinowska}
\date{}
\maketitle


\addtocounter{section}{-1}

\section{Introduction}\label{intro}
 
Soit $Y$ une vari\'et\'e complexe projective lisse de dimension $N+1$
munie d'un faisceau inversible tr\`es ample $\anneau (1)$~; soit 
$d\in\NN^*$~; soit ${\cal V}^{d}\subset\HH^0(Y,\anneau (d))$ l'ouvert
param\'etrant (\`a la multiplication par un scalaire pr\`es) les
hypersurfaces projectives lisses de $Y$ de classe $c_1(\anneau (d))$~;
soit $F\in{\cal V}^{d}$ d\'efinissant l'hypersurface $X_F$~; et soit
$j_F:X_F\to Y$ l'immersion ferm\'ee.  
Rappelons que la cohomologie \'evanescente de $X_F$ est l'espace 
$$\HH^{k}(X_F,\C)_{\ev}=
  \kker(g:\HH^{k}(X_F,\C)\rightarrow\HH^{k+2}(Y,\C)),$$ 
o\`u $k\in\{0,\dots,2N\}$ et o\`u $g$ d\'esigne l'application de 
Gysin, d\'eduite de l'application  
$j_{F*}:\HH_{2N-k}(X_F,\C)\rightarrow\HH_{2N-k}(Y,\C)$ 
par la dualit\'e de Poincar\'e.

Le groupe $\pi_1({\cal V}^{d},F)$ agit par monodromie sur
$\HH^{k}(X_F,\C)$. Cette action est d\'ecrite par la th\'eorie de
Lefschetz~:
\begin{itemize}
\item pour $k\neq N$ nous avons $\HH^{k}(X_F,\C)_{\ev}=0$ et la
repr\'esentation de monodromie de $\pi_1({\cal V}^{d},F)$ sur $\HH^{k}(X_F,\C)$
est triviale~;
\item les espaces $\HH^{N}(X_F,\C)_{\ev}$ et $j_F^*\HH^{N}(Y,\C)$ sont
en somme directe orthogonale pour la forme d'intersection~; la
repr\'esentation de monodromie de $\pi_1({\cal V}^{d},F)$ sur $\HH^{N}(X_F,\C)$
est somme directe de la repr\'esentation triviale sur 
$j_F^*\HH^{N}(Y,\C)$ et d'une repr\'esentation irr\'eductible sur
$\HH^{N}(X_F,\C)_{\ev}$.  
\end{itemize}

L'objectif de cet article est d'\'etendre ces r\'esultats au
cas d'une famille d'hypersurfaces lisses contenant un sous-sch\'ema
fix\'e. 

Plus pr\'ecis\'ement, soit $W\subset Y$ un sous-sch\'ema
ferm\'e de dimension $n$ d\'efini par un faisceau d'id\'eaux 
${\cal I}_W$. Nous supposons l'id\'eal homog\`ene
$I_W=\bigoplus_{i\in\NN}\HH^0(Y,{\cal I}_W\otimes\anneau(i))$ 
engendr\'e en degr\'e strictement inf\'erieur \`a un entier $e$. 
Soit ${\cal V}^{d}(W)\subset{\cal V}^{d}$ l'espace param\'etrant les
hypersurfaces lisses contenant $W$. Nous montrons d'abord

\begin{theoreme}\label{mono1}
Si $N<2n$, pour tout $d\geq e$, l'espace ${\cal V}^{d}(W)$ est vide.
\end{theoreme}

Supposons $N\geq 2n$. Bien s\^ur, l'espace ${\cal V}^{d}(W)$ peut encore
\^etre vide (c'est le cas si $W$ est trop singulier). Sinon, nous fixons
$F\in{\cal V}^{d}(W)$ et nous \'etudions l'action de
monodromie sur $\HH^{N}(X_F,\C)$ du groupe fondamental
$\pi_1({\cal V}^{d}(W),F)$. 

Soit $H(W)$ le sous-espace vectoriel de $\HH^{N}(X_F,\C)$ engendr\'e 
par $j_F^*\HH^{N}(Y,\C)$ et, lorsque $N=2n$, par les classes des 
composantes irr\'eductibles de $W$ de dimension~$n$. 
L'espace $H(W)$ est $\pi_1({\cal V}^{d}(W),F)$-invariant.
De plus, par le th\'eor\`eme d'indice de Hodge ({\it cf.}~\cite{weil}), la
forme d'intersection restreinte \`a $H(W)$ est non-d\'eg\'en\'er\'ee
(en effet, elle est non d\'eg\'en\'er\'ee sur 
$H(W)\cap\HH^{N}(X_F,\C)_{\ev}\subset\HH^{n,n}(X)$ d'une part, et sur 
$j_F^*\HH^{N}(Y,\C)$ d'autre part). Comme
l'action de $\pi_1({\cal V}^{d}(W),F)$ pr\'eserve la forme d'intersection, nous
avons l'\'egalit\'e de $\pi_1({\cal V}^{d}(W),F)$-modules 
$\HH^{N}(X_F,\C)=H(W)\oplus H(W)^{\perp}$. 
Le th\'eor\`eme principal est alors

\begin{theoreme}\label{mono}
Si $N\geq 2n$, il existe une constante $C\in\RR^*_+$ qui ne 
d\'epend que de $Y$, telle que pour tout $d\geq Ce$ et pour 
tout $F\in{\cal V}^{d}(W)$ la repr\'esentation de monodromie de
$\pi_1({\cal V}^{d}(W),F)$ sur $H(W)^{\perp}$ est irr\'eductible.
\end{theoreme}

Une hypoth\`ese sur $I_W$ est n\'ecessaire~: par exemple, si $d=1$ et si
$W$ est contenu dans un seul hyperplan strict de $Y$ d'\'equation $F$,
la conclusion du  th\'eor\`eme~\ref{mono} est fausse~: nous avons
${\cal V}^{d}(W)=\C^*\{F\}$ et la repr\'esentation de monodromie de
$\pi_1({\cal V}^{d},F)$ sur $\HH^{N}(X_F,\C)$ est triviale quelle que 
soit la dimension de $W$. En revanche, lorsque $I_W$ est engendr\'e en
degr\'e strictement inf\'erieur \`a $d$, nous montrons que le fait que $W$
soit contenu dans une hypersurface lisse de degr\'e $d$ est \'equivalent \`a
une condition portant sur le lieu singulier de $W$~; nous en d\'eduisons
que si $N<2n$, l'espace ${\cal V}^{d}(W)$ est vide (c'est-\`a dire le 
th\'eor\`eme~\ref{mono1}), et si $N\geq 2n$, 
l'espace ${\cal V}^{d}(W)$ est soit vide soit suffisamment gros (la
proposition~\ref{bertini} donne un \'enonc\'e pr\'ecis).

\medskip

Dans le cas o\`u $W$ est lisse, le th\'eor\`eme~\ref{mono} est vrai pour 
$C=1$ gr\^ace \`a l'argument tr\`es simple suivant, qui m'a \'et\'e communiqu\'e
par Voisin. 
Soit $\pi_W:\widetilde{Y_W}\to Y$ l'\'eclat\'e de $Y$ le long de $W$ et
soit $E_W$ le diviseur exceptionnel de $\widetilde{Y_W}$. 
Le transform\'e strict $\widetilde{X_W}$ d'une hypersurface lisse $X$
contenant $W$ est une hypersurface lisse de $\widetilde{Y_W}$,
isomorphe \`a l'\'eclat\'e de $X$ le long de $W$, de classe
$d\pi_W^*c_1(\anneau _{Y}(1))-c_1(E_W)$. Les sections du fibr\'e
de classe $d\pi_W^*c_1(\anneau _{Y}(1))-c_1(E_W)$ sont en bijection
avec $\HH^0(Y,{\cal I}_W\otimes\anneau(d))$, donc le fait que $I_W$ est engendr\'e en degr\'e strictement
inf\'erieur \`a $d$ implique que ce fibr\'e est tr\`es ample. L'\'enonc\'e
r\'esulte alors par le th\'eor\`eme de Lefschetz difficile de
l'irr\'eductibilit\'e de la repr\'esentation de monodromie de
$\pi_1({\cal V}^{d}(W))$ sur la cohomologie \'evanescente de $\widetilde{X_W}$.  

Malheureusement cette preuve ne s'\'etend
pas au cas singulier. Supposons en effet que $W$ poss\`ede un point
singulier $w$ qui n'est pas un point double ordinaire. Pour un pinceau
g\'en\'erique $\pi:{\cal X}\to L$, o\`u 
$L\subset\PP\HH^0(Y,{\cal I}_W\otimes\anneau(d))$, il peut exister
un point $l\in L$ tel que le sch\'ema ${\cal X}$ est singulier au point
$(w,l)$ et $(w,l)$ n'est pas un point double ordinaire. La monodromie
locale de $\pi_1({\cal V}^{d}\cap L)$ au voisinage de $l$ ne v\'erifie alors pas
n\'ecessairement la formule de Picard-Lefschetz. Contrairement au cas o\`u
$W$ est lisse, il n'existe donc pas en g\'en\'eral d'\'eclatement 
$\pi_{\cal X}:\widetilde{\cal X}\to {\cal X}$ tel que la famille
$\pi\circ\pi_{\cal X}:\widetilde{\cal X}\to L$ soit un pinceau de
Lefschetz. Nous pensons cependant que le th\'eor\`eme~\ref{mono} reste
vrai avec $C=1$ lorsque $W$ est singulier.

\medskip

La motivation principale du th\'eor\`eme~\ref{mono} appara\^it dans
son application \`a la th\'eorie de Hodge. Rappelons que la
conjecture de Hodge implique que pour tout $F\in{\cal V}^{d}$,
l'espace des classes de Hodge \'evanescentes 
${\rm Hdg(X_F)}_{\ev}=\HH^{n,n}(X_F)\cap\HH^{2n}(X,\QQ)_{\ev}$
est engendr\'e par la projection sur la cohomologie 
\'evanescente des classes 
des sous-vari\'et\'es alg\'ebriques de $X_F$.
Le lieu $NL_d=\{ F\in{\cal V}^{d}\ |\ {\rm Hdg(X_F)}_{\ev}\neq 0\}$
o\`u cet \'enonc\'e est non vide
s'appelle le {\em lieu de Noether-Lefschetz}. Si $d$ est assez grand 
c'est une r\'eunion d\'enombrable de sous-vari\'et\'es alg\'ebriques 
strictes de ${\cal V}^{d}$. Nous montrons dans un autre article que
pour $d\gg 0$ une hypersurface g\'en\'erique appartenant \`a une 
composante de $NL_d$ de dimension suffisamment grande
v\'erifie la conjecture de Hodge ci-dessus.
Le th\'eor\`eme~\ref{mono} est un ingr\'edient essentiel de la preuve.

\medskip

La preuve du th\'eor\`eme~\ref{mono} repose sur la construction d'une
filtration de la (co)homologie de $X_F$ au voisinage de la
d\'eg\'en\'erescence de $X_F$ en la r\'eunion de deux hypersurfaces lisses se
coupant transversalement  
(proposition~\ref{decompositionpremiere}). Notre d\'emarche s'inspire
de~\cite{hodgevariationnel} et~\cite{lopez} qui \'etudient le lieu de
Noether-Lefschetz des surfaces de $\PP^3_{\C}$ par une m\'ethode
similaire. Contrairement aux filtrations \'etudi\'ees habituellement, et
notamment \`a la filtration par le poids \'etudi\'ee par Clemens~\cite{clem},
Schmid~\cite{sm} et Steenbrink~\cite{ste},
notre filtration introduit une dissym\'etrie entre les hypersurfaces
constituant la  fibre singuli\`ere. En particulier, la cohomologie
\'evanescente de l'une des hypersurfaces n'est pas un sous-espace
du gradu\'e associ\'e \`a notre filtration (mais en est seulement un
sous-quotient)~; cette propri\'et\'e est cruciale pour la 
d\'emonstration du th\'eor\`eme~\ref{mono}. 

La suite de la preuve du th\'eor\`eme~\ref{mono} est la suivante. Nous
fixons $i\in\NN$ tel que $e\leq i<d$ et nous \'etudions la (co)homologie
de $X_F$ au voisinage de sa d\'eg\'en\'erescence en $X_A\cup X_K$ avec 
$A\in{\cal V}^i(W)$ et $K\in{\cal V}^{d-i}$. Nous exhibons une action
naturelle de monodromie de $\pi_1({\cal V}^{d-i},K)$ sur tous les
gradu\'es associ\'es \`a la filtration de la (co)homologie de $X_F$ que nous
avons construite. Cette action se factorise par l'action de monodromie
de $\pi_1({\cal V}^{d}(W),F)$. Sur tous les gradu\'es sauf un les deux
repr\'esentations ont m\^emes orbites et nous pouvons minorer les
dimensions minimales des sous-repr\'esentations dans ces gradu\'es par
les dimensions des espaces de (co)homologie \'evanescente de $X_K$ et de
$X_A\cap X_K$, sur lesquelles l'action de $\pi_1({\cal V}^{d-i},K)$ est
irr\'eductible. Pour d\'ecrire l'action de $\pi_1({\cal V}^{d}(W),F)$ sur le
dernier gradu\'e 
nous proc\'edons par r\'ecurrence sur $N$. Nous montrons ainsi que la
dimension de la $\pi_1({\cal V}^{d}(W),F)$-repr\'esentation de monodromie
engendr\'ee par une classe de $H(W)^{\perp}$ est minor\'ee par l'homologie
\'evanescente d'une intersection compl\`ete lisse dans $X_F$
d'hypersurfaces de degr\'es $i$ et $d-i$ (proposition~\ref{aux}). Nous
concluons par des estimations num\'eriques de la dimension de ces espaces
(proposition~\ref{estimation}), en faisant varier l'entier $i$.

\medskip

\noindent
{\bf Remerciements.}
Je remercie Morihiko Saito pour sa lecture tres attentive de mon article.
Il m'a signal\'e deux impr\'ecisions, une dans la preuve du 
lemme~\ref{contraction} et une (plus s\'erieuse) dans la 
section~\ref{preuveaux}, et m'a indiqu\'e la r\'ef\'erence~\cite{ka} 
qui simplifie la preuve de la proposition~\ref{bertini}.
Cette version tient compte de ses remarques.
D'autre part, il a trouve une variante de la preuve du th\'eor\`eme~\ref{mono}
qui ameliore notablement la borne $C$. Ce travail sera l'objet d'un article 
commun, ecrit dans un langage beaucoup plus sophistiqu\'e. 


\section{Sch\'emas contenus dans une
hypersurface lisse de grand degr\'e~; preuve du th\'eor\`eme~\ref{mono1}}

\subsection{Caract\'erisation intrins\`eque et propri\'et\'es}

Dans cette partie nous caract\'erisons les   
sous-sch\'emas ferm\'es $W\subset Y$ contenus dans une hypersurface
lisse suffisamment ample de $Y$ par une condition portant uniquement
sur le lieu singulier de $W$. Nous \'etablissons ensuite un \'enonc\'e de
type {\it Bertini} pour la famille d'hypersurfaces contenant $W$ et nous
montrons que dans une intersection compl\`ete g\'en\'erique d'un nombre
maximal d'hypersurfaces suffisamment amples contenant $W$ la
{\it liaison} de $W$ est connexe et lisse. L'arbitre nous signale que
ces r\'esultats se trouvent partiellement d\'emontr\'es dans \cite{ka}, 
et nous l'en remercions.

Plus pr\'ecis\'ement, posons $I_W^i=\HH^0(Y,{\cal I}_W\otimes\anneau_Y(i))$ 
et notons $I_W=\bigoplus_{i\in\NN}I_W^i$ l'id\'eal homog\`ene
d\'efinissant $W$ dans $Y$~; soit $e$ le plus petit entier tel que
$I_W$ est engendr\'e en degr\'e strictement inf\'erieur \`a $e$.    
Soit enfin $d$ un entier. Nous nous demandons \`a quelle condition sur $W$ il
existe une hypersurface lisse $X\subset Y$ de classe $c_1(\anneau(d))$
contenant $W$. 

Si $d<e$ le probl\`eme n'a pas de solution simple, comme l'illustre le
cas $Y=\PP^{N+1}_{\C}$, $d=1$~: nous imposons seulement que $W$ soit
r\'ealisable comme un sous-sch\'ema de $\PP^{N}_{\C}$, ce qui autorise
des singularit\'es arbitraires, sauf la condition \'evidente que 
pour tout $x\in W$ l'espace tangent \`a $W$ en $x$ doit \^etre de 
dimension inf\'erieure ou \'egale \`a $N$. Notons que dans ce cas le 
sch\'ema $W$ peut \^etre contenu dans une seule hypersurface lisse.

Si $d\geq e-1$, c'est-\`a-dire si $I_W$ est engendr\'e en degr\'e 
inf\'erieur ou \'egal \`a $d$, Kleinman et Altman montrent qu'il
existe une condition suffisante simple portant sur la dimension des  
espaces tangents \`a $W$ pour qu'une hypersurface g\'en\'erique de classe 
$c_1(\anneau(d))$ contenant $W$ soit lisse. 

Si $d\geq e$, nous montrons ci-dessous que cette condition devient
n\'ecessaire. Dans ce cas le fait que $W$ soit contenu dans une
hypersurface lisse de classe $c_1(\anneau(d))$ est donc enti\`erement
d\'etermin\'e par la dimension des espaces tangents \`a $W$,
contrairement au cas $d<e$.   

\medskip

Plus pr\'ecis\'ement, pour tout $c\in \{0,\dots, N-n+1\}$ posons 
$$W_c=\{x\in W\ |\ \codim(T_xW,T_xY)=c\},$$
o\`u $T_xW$ d\'esigne l'espace tangent de Zariski \`a $W$ en $x$~: si
${\cal I}_{W,x}$ d\'esigne l'id\'eal de l'anneau local 
${\cal O}_{Y,x}$ et $d$ la diff\'erentielle, alors
$T_xW=\bigcap_{q\in{\cal I}_{W,x}}\kker dq_x$.
L'espace $W_c$ est un ensemble alg\'ebrique localement ferm\'e.

Consid\'erons les assertions suivantes, o\`u $\delta$ est un entier. 

{\it {\renewcommand{\labelenumi}{$($\alph{enumi}$)$}
\begin{enumerate}
\item 
Il existe un entier $d\geq \delta$ et une hypersurface lisse de $Y$ de
classe $c_1(\anneau(d))$ contenant $W$.
\item 
Pour tout $c\in \{0,\dots, N-n+1\}$ nous avons $\dim W_c<c$.
\item
Nous avons $N\geq 2n$. De plus, pour tout $r\in\{1,\dots,N-n+1\}$, pour 
tout $r$-uplet $(e_1,\dots,e_r)$ d'entiers sup\'erieurs ou
\'egaux \`a $\delta$ et pour $(P_1,\dots,P_r)\in
\prod_{i=1}^{r}I_W^{e_i}$ g\'en\'erique, le $r$-uplet
$(P_1,\dots,P_r)$ d\'efinit une intersection compl\`ete et le lieu
singulier du sch\'ema $Z(P_1,\dots,P_{r})$ est support\'e par $W$ et est
de dimension inf\'erieure ou \'egale \`a $r-2$ (par convention, le lieu
singulier d'une vari\'et\'e lisse est de dimension $-1$).
\item 
Pour tout $r\in\{1,\dots,N-n+1\}$, $r\neq N+1$, pour tout
$r$-uplet $(e_1,\dots,e_r)$ d'entiers sup\'erieurs ou \'egaux \`a 
$\delta$ et pour $(P_1,\dots,P_r)\in \prod_{i=1}^{r}I_W^{e_i}$ 
g\'en\'erique le sch\'ema $Z(P_1,\dots ,P_{r})\setminus W$ est 
connexe et lisse. 
\end{enumerate} }}

L'implication $(c)\Rightarrow(a)$ est \'evidente.

D'apr\`es \cite{ka}, nous avons 
$(b)\Rightarrow (c)$ pour $\delta\geq e-1$ et $(b)\Rightarrow (d)$ 
pour $\delta\geq e$ et $r\neq N-n+1$.

\begin{prop} \label{bertini} 
Pour $\delta\geq e$, les assertions 
$(a)$, $(b)$ et $(c)$ sont \'equivalentes et elles impliquent
l'assertion $(d)$.
\end{prop}

L'implication $(a)\Rightarrow (c)$ donne le th\'eor\`eme~\ref{mono1}. 

\subsection{Preuve de la proposition~\ref{bertini}}
Il suffit de montrer  $(a)\Rightarrow (b)$ et $(b)\Rightarrow (d)$
pour $r=N-n+1$.

\medskip

\subsubsection{Preuve de $(a)\Rightarrow (b)$}

Soit $x\in W_c$, soit $L\in\HH^0(Y,\anneau(1))$ une forme lin\'eaire
ne s'annulant pas 
en $x$ et soit $U$ un voisinage de Zariski ouvert de $x$ dans $W_c$
tel que $L$ ne s'annule pas sur $U$. En choisissant $x\in W_c$ dans
une composante irr\'eductible de dimension maximale de $W_c$ nous
pouvons supposer $U$ de m\^eme dimension que $W_c$. Posons
$$E=\coker(TW_{|W_c}\to TY_{|W_c})_{|U}^{\vee}.$$ 
C'est un fibr\'e vectoriel de rang $c$ sur $U$.

L'espace $E_x$ s'identifie canoniquement \`a $(T_xY/T_xW)^{\vee}$~; il 
est donc engendr\'e par les formes diff\'erentielles $d_xq$ o\`u $q$
d\'ecrit l'id\'eal ${\cal I}_{W,x}$. Cet id\'eal est engendr\'e en
tant que ${\cal O}_{Y,x}$-module par les \'el\'ements de la forme
$\frac{Q}{L^{\deg\,Q}}$, o\`u $Q\in I_W$ est un \'el\'ement
homog\`ene. L'espace $\HH^0(U,E)$ est donc engendr\'e par l'image des 
applications $d^{L,i}:I_W^{i}\to \HH^0(U,E)$, 
$Q\mapsto d\left(\frac{Q}{L^{i}}\right)$ o\`u $i$ d\'ecrit $\NN$. 
Remarquons que $d^{L,i}Q$ s'annule en $x$ si et seulement si 
l'hypersurface $X_Q$ est singuli\`ere en $x$.

Comme $I_W$ est engendr\'e en degr\'e inf\'erieur ou \'egal \`a
$\delta-1$, donc \`a $d-1$, l'image de  $d^{L,d-1}$ dans $\HH^0(U,E)$
engendre le fibr\'e $E$. Donc l'image de $d^{L,d}$ engendre les 1-jets
de $E$, {\it i.e.} l'image de $d^{L,d}$ engendre $E$ et pour tout
$x\in E$ les diff\'erentielles des \'el\'ements de l'image de
$d^{L,d}$ qui s'annulent  en $x$ engendrent 
$E_x\otimes T_xU^{\vee}$. En effet, si $L'\in\HH^0(Y,\anneau(1))$ est
une forme lin\'eaire s'annulant en $x$ et si $Q\in I_W^{d-1}$, alors
la section $d^{L,d}(QL')$ s'annule en $x$ et  
$d_x(d^{L,d}(QL'))=d^{L,d-1}Q\otimes d_x\left(\frac{L'}{L}\right)$~;
l'\'enonc\'e r\'esulte alors de ce que les $d^{L,d-1}Q$ engendrent $E$
et les $d_x\left(\frac{L'}{L}\right)$ engendrent $T_xU^{\vee}$.

Quitte \`a restreindre $U$ nous pouvons le supposer trivial~; l'espace 
$\HH^0(U,E)$ s'identifie alors aux fonctions alg\'ebriques de $U$ dans
$\C^c$. Par l'hypoth\`ese $(a)$, il existe $Q\in I_W^d$ tel que $X_Q$ 
est lisse. Comme la lissit\'e est une propri\'et\'e ouverte, pour 
$Q\in I_W^d$ g\'en\'erique $X_Q$ est lisse~; la fonction
$d^{L,d}Q$ ne s'annule nulle part sur U. L'assertion $(b)$ r\'esulte
donc par contrapos\'ee du lemme suivant appliqu\'e \`a 
$\Sigma=d^{L,d}(I_W^{d})$.
\cucu

\begin{lemme}
Soit $U$ une vari\'et\'e lisse de dimension sup\'erieure ou 
\'egale \`a $c$ et $\Sigma$ un sous-espace de l'espace des
fonctions alg\'ebriques de $U$ dans $\C^c$ qui engendre 
les 1-jets. Alors pour $\sigma\in\Sigma$ g\'en\'erique, le lieu
$\sigma^{-1}(0)$ est non vide.
\end{lemme}

\preuve  
Comme $\Sigma$ engendre les 1-jets et que $\dim X \geq c$, il existe
$(\sigma,x)\in\Sigma\times U$ tels que $\sigma(x)=0$ et la
diff\'erentielle  $d_x\sigma:T_x U\rightarrow \C^c$ est surjective.  

Notons ${\rm ev}\colon\Sigma\times U\to\C^c$ l'application 
$(\tau,y)\mapsto\tau(y)$ et posons $\Theta={\rm ev}^{-1}(0)$.
Nous avons $(\sigma,x)\in\Theta$.
Consid\'erons alors le diagramme commutatif dont les lignes sont exactes
$$\xymatrix{0\rto&
       T_xU\dto^{d_x\sigma}\rto^{d_x(\sigma,\cdot)\quad}&
       T_{\sigma,x}\Sigma\times U\dto^{d_{\sigma,x}{\rm ev}}
                   \rto^{\quad d_{\sigma,x}\pi}&
       T_{\sigma}\Sigma\dto\rto&0\\
            0\rto&\C^c\rto&\C^c\rto&0\rto&0\\ }$$
o\`u $\pi\colon\Sigma\times U\to\Sigma$ d\'esigne la projection. 
Par chasse au diagramme, la surjectivit\'e de $d_x\sigma$ implique
que la restriction de $d_{\sigma,x}\pi$ \`a $\kker\,d_{\sigma,x}{\rm ev}=
T_{\sigma,x}\Theta$ est surjective.
\cucu

\subsubsection{Preuve de $(b)\Rightarrow (d)$}
D'apr\`es \cite{ka}, pour
$(P_1,\dots,P_{N-n})\in\Pi_{i=1}^{N-n}I_W^{e_i}$ 
g\'en\'erique le sch\'ema $Z(P_1,\dots,P_{N-n})\setminus W$ est
connexe et lisse. Fixons $(P_1,\dots,P_{N-n})$ comme ci-dessus, 
choisissons $P_{N-n+1}\in I_W^{e_{N-n+1}}$, posons 
$Y'=Z(P_1,\dots,P_{N-n})$ et $X'=Y'\cap Z(P_{N-n+1})$. Nous devons    
montrer que pour $P_{N-n+1}\in I_W^{e_{N-n+1}}$ g\'en\'erique le
sch\'ema $X'\setminus W$ est connexe et lisse. 

Notons $\pi':\widetilde{Y'}\to Y'$ l'\'eclat\'e de 
$Y'$ le long de $W$, $E'$ le diviseur
exceptionnel et $\widetilde{X'}$ le transform\'e strict de $X'$.
Ainsi $\widetilde{X'}$ est une hypersurface de
$\widetilde{Y'}$, isomorphe \`a l'\'eclat\'e de $X'$ le long de
$W$, section du fibr\'e ${\cal L}'$ sur $\widetilde{Y'}$ de classe 
$e_{N-n+1}\pi^{'*}c_1(\anneau _{Y}(1))-c_1(E')$. Comme les sections de
${\cal L}'$ sont en bijection avec $I_W^{e_{N-n+1}}$ et que $I_W$ est 
engendr\'e en degr\'e $e-1<\delta\leq e_{N-n+1}$,  
le fibr\'e ${\cal L}'$ est tr\`es ample. D'autre part, le sch\'ema 
$\widetilde{Y'}\setminus E'\simeq Y'\setminus W$ est
connexe et lisse, de dimension $n+1\geq 2$ (puisque ${N-n+1}\neq N+1$); le 
th\'eor\`eme de Bertini classique ({\it cf}. par exemple \cite{ka})
implique alors que pour $P_{N-n+1}\in I_W^{e_{N-n+1}}$ g\'en\'erique le sch\'ema 
$X'\setminus (E'\cap X')\simeq X'\setminus W$ est 
irr\'eductible et lisse, donc connexe et lisse.
\cucu

\section{Les r\'esultats de monodromie}

\subsection{Notations et rappels sur la th\'eorie de Lefschetz}\label{m1}

Dans cette section nous rappelons des r\'esultats dus essentiellement \`a
Lefschetz (sauf le th\'eor\`eme de Lefschetz difficile) et expliqu\'es en
d\'etail dans~\cite{lamotke}. Nous les utiliserons
librement dans les parties 2 et 3.  

Le th\'eor\`eme~\ref{mono} affirme essentiellement l'irr\'eductibilit\'e 
d'une action de monodromie~; nous l'\'enon\c cons pour l'homologie 
complexe, qui donne le r\'esultat le plus fort; nous nous pla\c cons donc 
dans ce cadre~: sauf mention explicite du contraire, tous les espaces
d'homologie consid\'er\'es dans les parties 2 et 3 sont \`a coefficients
complexes.  

\subsubsection{Homologie primitive, \'evanescente, relative} 

\noindent{\it Homologie primitive.\ ---\ }
Soit $N$ un entier strictement positif et soit
$Y$ une vari\'et\'e projective lisse de dimension 
$N+1$ munie d'un faisceau inversible tr\`es ample ${\cal L}$. Nous 
appelons homologie primitive de $Y$ les espaces
\begin{eqnarray*}
\HH_i(Y)^{\prim}&=&\coker(c_1({\cal L}):
\HH_{i+2}(Y)\rightarrow\HH_i(Y))\quad {\rm et}\\
\HH_i(Y)_{\prim}&=&\kker(c_1({\cal L}):
\HH_{i}(Y)\rightarrow\HH_{i-2}(Y)),
\end{eqnarray*}
et nous posons $\hh_{i}(Y)^{\prim}=\dim\HH_i(Y)^{\prim}$. 
Nous avons $\HH_{2N+2-i}(Y)_{\prim}=\HH_{i}(Y)^{\prim}=0$ 
pour $i\in\{0,\dots,N\}$.

\medskip

\noindent{\it Homologie \'evanescente.\ ---\ }
Soit $X$ une hypersurface projective lisse de $Y$ de classe 
$c_1({\cal L})$ et soit $j:X\to Y$ l'immersion ferm\'ee.
Nous appellons homologie \'evanescente les espaces
\begin{eqnarray*}
\HH_i(X)_{\ev}&=&
\kker(j_*:\HH_i(X)\rightarrow\HH_{i}(Y))\quad {\rm et}\\
\HH_i(X)^{\ev}&=&
\coker(j^*:\HH_{i+2}(Y)\rightarrow\HH_i(X)).
\end{eqnarray*}

L'homologie \'evanescente de $X$ d\'epend \`a priori du plongement
$j:X\to Y$. Cependant, si $X$ d\'efinit un diviseur tr\`es ample dans 
deux vari\'et\'es lisses $Y$ et $Y'$, et s'il existe une
vari\'et\'e lisse $Z$ telle que $Y\subset Z$ et $Y'\subset Z$ sont des
produits d'intersection de diviseurs tr\`es amples, alors 
les homologies \'evanescentes de $X$ pour les plongements 
$X\to Y$ et $X\to Y'$ co\"{\i}ncident. Dans cet article, nous nous
trouvons toujours dans cette situation. 

Nous avons $\HH_i(X)_{\ev}=\HH_i(X)^{\ev}=0$ pour
$i\in\{0,\dots,2N\}\setminus \{N\}$. La forme d'intersection 
induit un isomorphisme canonique
$\HH_N(X)_{\ev}\simeq\HH_N(X)^{\ev}$~;
nous notons $\hh_N(X)_{\ev}$ la dimension de cet espace.  

\medskip

\noindent{\it Homologie relative et homologie du compl\'ementaire.\ ---\ }
Notons ${\cal V}^i\subset\HH^0(Y,{\cal L})$
l'ouvert param\'etrant (\`a la multiplication par un scalaire pr\`es)
les hypersurfaces lisses de classe $c_1({\cal L})$.   
Pour tout $i\in\{0,\dots,2N\}$ nous avons les suites exactes de
$\pi_1({\cal L})$-modules 
\begin{eqnarray}
&0\longrightarrow \HH_{i+1}(Y)^{\prim}\stackrel{\rel}{\longrightarrow}
\HH_{i+1}(Y,X)\stackrel{\bord}{\longrightarrow} 
\HH_{i}(X)_{\ev} \longrightarrow 0, \label{rela}&\\
&0\longrightarrow \HH_{i}(X)^{\ev}
\stackrel{\tube}{\longrightarrow} \HH_{i+1}(Y\setminus X)
\stackrel{\ouv}{\longrightarrow}\HH_{i+1}(Y)_{\prim}
\longrightarrow 0.&
\label{ouvert} 
\end{eqnarray}
Ces suites exactes sont duales
pour les formes d'intersection sur $X$ et sur $Y$~; en 
particulier, nous avons des isomorphismes canoniques  
$\HH_{i+1}(Y,X)\simeq\HH_{2N-i+1}(Y\setminus X)^{\vee}.$

\subsubsection{Monodromie}\label{rappelmonodromie}
 
Nous introduisons les notations suivantes. Si $B_*$ d\'esigne un
ouvert de $Y$ nous notons $\Gamma_*$ son compl\'ementaire,
$\overline{B}_*$ son adh\'erence et $\Sigma_*$ le bord de
$\overline{B}_*$. Pour tout sous-espace $Z$ de $Y$ nous posons 
$B_*^{Z}=Z\cap B_*$, $\overline{B}_*^{Z}=Z\cap\overline{B}_*$,
$\Gamma_*^{Z}=Z\cap\Gamma_*$ et $\Sigma_*^{Z}=Z\cap\Sigma_*$. Nous
notons $\beta_{B,Z}:\HH_{\bullet}(B_*^Z)\to\HH_{\bullet}(Z)$ le 
morphisme induit par l'immersion ferm\'ee $B_*\to Z$ et 
$\gamma_{B,Z}:\HH_{\bullet}(Z)\to\HH_{\bullet}(\overline{B}_*^Z,\Sigma_*^Z)$ la
compos\'ee du morphisme 
$\rel:\HH_{\bullet}(Z)\to\HH_{\bullet}(Z,\Gamma_*^Z)$ avec
l'isomorphisme d'excision  
$\HH_{\bullet}(Z,\Gamma_*^Z)\to\HH_{\bullet}(\overline{B}_*^Z,\Sigma_*^Z)$.

\medskip

Soit $D\subset\PP\HH^0(Y,{\cal L})$ un pinceau de Lefschetz 
d'hypersurfaces de $Y$. Notons $\{P_1,\dots, P_r\}\subset D$ ses 
points critiques~: pour tout $i\in\{1,\dots, r\}$ 
l'hypersurface associ\'ee \`a $P_i$ poss\`ede un unique point double 
ordinaire not\'e $x_i$. Posond $D^*=D\setminus\{P_i,\dots, P_r\}$

\medskip

Nous fixons $i\in\{1,\dots, r\}$. Soit $B_i$ un voisinage de $x_i$
dans $Y$  hom\'eomorphe \`a une boule ouverte.
D'apr\`es la th\'eorie de Morse, pour $B_i$ assez petit, il existe 
un voisinage $D_i$ de $P_i$ dans $D$ hom\'eomorphe \`a un disque 
ferm\'e et ne rencontrant pas les $P_j,\ j\not=i$ tel que 
la famille de vari\'et\'es \`a bord $(\Gamma_i^{X},\Sigma_i^{X})$
est topologiquement triviale pour $X$ d\'ecrivant la famille
d'hypersurfaces param\'etr\'ee par $D_i$, et que la famille $B_i^{X}$
est topologiquement localement triviale pour $X$ d\'ecrivant la
famille d'hypersurfaces param\'etr\'ee par $D_i\setminus\{P_i\}$. 

Soit $F_i$ un point de $D_i\setminus\{P_i\}$ d\'efinissant une  
hypersurface $X_i\subset Y$~: la vari\'et\'e $B_i^{X_i}$ est alors 
hom\'eomorphe \`a un sous-fibr\'e en boules ouvertes de dimension $N$ 
du  fibr\'e tangent \`a la sph\`ere r\'eelle ${\cal S}^N$ et  
la vari\'et\'e $\Sigma_i^{X_i}$ est hom\'eomorphe \`a un sous-fibr\'e en 
sph\`eres r\'eelles de dimension $N-1$ du fibr\'e tangent \`a 
${\cal S}^N$. L'espace $\HH_N(B_i^{X_i})$ est donc isomorphe \`a $\C$ 
et est engendr\'e par la classe de l'image de la section nulle, que 
nous notons $\delta_i^{B_i}$. 
L'espace $\HH_N(\overline{B}_i^{X_i},\Sigma_i^{X_i})$ est le dual de 
$\HH_N(B_i^{X_i})$ pour la forme d'intersection, engendr\'e par la
classe d'une boule-fibre, qui donc coupe transversalement la section
nulle. L'application $\gamma_{B_i,X_i}$ est la transpos\'ee de
$\beta_{B_i,X_i}$. 

Nous posons $\delta_i=\beta_{B_i,X_i}(\delta_i^{B_i})$. Soit $W_i \subset
D_i\setminus P_i$ un lacet d'origine $F_i$ faisant le tour de $P_i$ et
soit $g_{i}\in\pi_1(D_i\setminus P_i,F_i)$ sa classe d'homotopie. 
La formule de Picard-Lefschetz affirme que pour tout
$\lambda\in\HH_{N}(X_i)$ nous avons  
$$g_i(\lambda)=
  \lambda+\varepsilon\langle\lambda|\delta_i\rangle\delta_i,
  \quad {\rm avec}\quad \varepsilon=1\ {\rm ou}\ -1.$$

\medskip

Soit $F$ un point de $D^*$ d\'efinissant 
une hypersurface $X$~; pour tout $i\in\{1,\dots, r\}$ soit
$L_i\subset D^*$ un chemin reliant 
$F$ \`a $F_i$, induisant par transport plat un isomorphisme
$\HH_{N}(X)\simeq\HH_{N}(X_i)$~; nous appelons {\it cycle \'evanecsent}
associ\'e \`a $P_i$ et nous notons abusivement $\delta_i\in\HH_{N}(X)$ la
pr\'eimage de $\delta_i\in\HH_{N}(X_i)$ et
$g_i\in\pi_1(D^*,F)$ la classe d'homotopie 
du lacet $L_i^{-1}W_iL_i$. La formule de Picard-Lefschetz implique
alors que pour tout $\lambda\in\HH_{N}(X)$ et pour tout 
$i\in\{1,\dots, r\}$ nous avons 
$g_i(\lambda)=
 \lambda+\varepsilon\langle\lambda|\delta_i\rangle\delta_i$,
avec $\varepsilon=1$ ou $-1$.
La th\'eorie de Lefschetz affirme d'autre part
$\delta_i\in\HH_{N}(X_i)_{\ev}$ et que l'homologie \'evanescente  
$\HH_{N}(X)_{\ev}$ est engendr\'ee par les cycles 
$\delta_i$, $i\in\{1,\dots, r\}$ et que ces cycles sont conjugu\'es
sous l'action de monodromie de $\pi_1(D^*,F)$.
Ceci implique que l'action de monodromie de 
$\pi_1(D^*,F)$ sur $\HH_{N}(X)_{\ev}$ 
est irr\'eductible. 

\subsection{Quelques r\'esultats de monodromie globale}

\subsubsection{Monodromie du compl\'ementaire d'une hypersurface}
\label{m2}

Nous adoptons les notation de la section~\ref{rappelmonodromie}.
Le groupe $\pi_1({\cal V},F)$ 
agit par mo\-no\-dro\-mie sur les suites exactes (\ref{rela}) et
(\ref{ouvert}). Son action est triviale sur  
$\HH_{N+1}(Y)_{\prim}$ et sur $\HH_{N+1}(Y)^{\prim}$ et elle est
irr\'eductible sur $\HH_{N}(X)_{\ev}$ et sur $\HH_{N}(X)^{\ev}$.
La suite exacte (\ref{rela}) montre que les seules classes
$\pi_1({\cal V},F)$-invariantes dans $\HH_{N+1}(Y,X)$ sont celles 
appartenant \`a $\rel(\HH_{N+1}(Y)^{\prim})$, et que toutes les autres 
classes engendrent une $\pi_1({\cal V},F)$-repr\'esentation de dimension 
sup\'erieure ou \'egale \`a $\hh_{N}(X)_{\ev}$. L'\'enonc\'e
analogue pour la suite exacte (\ref{ouvert}) est l'objet de la
proposition suivante. 

\begin{prop}\label{monodromie-sur-complementaire}
Une classe non nulle $\Lambda\in\HH_{N+1}(Y\setminus X)$ engendre 
la $\pi_1({\cal V},F)$-repr\'e\-sen\-ta\-tion
$\C\Lambda+\tube\,\HH_{N}(X)^{\ev}$.
\end{prop}

\preuve 
Nous pouvons supposer $F\in D^*$~; il suffit alors de montrer le
r\'esultat pour $\pi_1(D^*,F)$ au lieu de $\pi_1({\cal V},F)$
Nous fixons $i\in\{1,\dots, r\}$. Consid\'erons la la suite
exacte longue d'homologie relative du couple $(B_i^{X_i},B_i)$~: comme 
$\HH_{N+1}(B_i)=\HH_{N}(B_i)=0$, nous avons un isomorphisme
$\bord_{B_i}:\HH_{N+1}(B_i,B_i^{X_i})\to\HH_{N}(B_i^{X_i})$.
Nous notons $\Delta_i^{B_i}$ la pr\'eimage de $\delta_i^{B_i}$ par
$\bord_{B_i}$. 

Soit $\beta_{B_i,Y,X_i}:\HH_{N+1}(B_i,B_i^{X_i})\to\HH_{N+1}(Y,X_i)$
le morphisme induit par les inclusions $B_i\to Y$ et $B_i^{X_i}\to X_i$
(le morphisme $\beta_{B_i,Y,X_i}$ est le transpos\'e de 
$\gamma_{B_i,Y\setminus X}:\HH_{N+1}(Y\setminus X_i)\to
\HH_{N+1}(\overline{B}_i^{Y\setminus X_i},\Sigma_i^{Y\setminus X_i})$
pour la forme d'intersection). Nous notons $\Delta_i$ l'image de
$\Delta_i^{B_i}$ par $\beta_{B_i,Y,X_i}$. Nous avons 
$\bord\,\Delta_i=\delta_i$.  

\begin{lemme}[Formule de Picard-Lefschetz pour le compl\'ementaire]
\label{picard-lefschetz}
Pour tout $\Lambda\in\HH_{N+1}(Y\setminus X_i)$ nous avons  
$$g_i(\Lambda)=
\Lambda+\varepsilon\langle\Lambda|\Delta_i\rangle\tube\,(\delta_i),
  \quad {\rm avec}\quad \varepsilon=1\ {\rm ou}\ -1.$$ 
\end{lemme}

\preuve
Nous distiguons les cas $\delta_i=0$ et $\delta_i\neq 0$.

Si $\delta_i=0$ alors ${\rm tube}\,\delta_i=0$~; nous devons donc
montrer $g_i={\rm id}$. Consid\'erons la suite exacte longue
d'homologie relative du couple $(B_i^{Y\setminus X_i},B_i)$~: 
comme $\HH_{N+2}(B_i)=\HH_{N+1}(B_i)=0$ et gr\^ace \`a l'isomorphisme de
Thom $\HH_{N+2}(B_i,B_i^{Y\setminus X_i})\simeq\HH_{N}(B_i^{X_i})$, nous
avons un isomorphisme
$\tube_{B_i}:\HH_{N}(B_i^{X_i})\to\HH_{N+1}(B_i^{Y\setminus X_i})$. Comme
l'espace $\HH_{N}(B_i^{X_i})$ est engendr\'e par $\delta_i^{B_i}$, l'espace
$\HH_{N+1}(B_i^{Y\setminus X_i})$ est donc engendr\'e par 
$\tube_{B_i}(\delta_i^{B_i})$. Or nous avons
$\beta_{B_i,Y\setminus X_i}(\tube_{B_i}(\delta_i^{B_i}))=
 \tube(\beta_{B_i,X_i}(\delta_i^{B_i}))=\tube(\delta_i)=0$~; le morphisme
$\beta_{B_i,Y\setminus X_i}:\HH_{N+1}(B_i^{Y\setminus X_i})
\to\HH_{N+1}(Y\setminus X_i)$ est donc nul. Le morphisme 
$\rel:\HH_{N+1}(Y\setminus X_i)\to
 \HH_{N+1}(Y\setminus X_i,B_i^{Y\setminus X_i})$ est donc injectif~;
en le composant avec l'isomorphisme d'excision
$\HH_{N+1}(Y\setminus X_i,B_i^{Y\setminus X_i})\to
 \HH_{N+1}(\Gamma_i^{Y\setminus X_i},\Sigma_i^{Y\setminus X_i})$,
nous obtenons un morphisme injectif $\HH_{N+1}(Y\setminus X_i)\to
\HH_{N+1}(\Gamma_i^{Y\setminus X_i},\Sigma_i^{Y\setminus X_i})$.
Le lemme r\'esulte alors de ce que la famille de vari\'et\'es \`a bord
$(\Gamma_i^{Y\setminus X_t},\Sigma_i^{Y\setminus X_t})$  
est topologiquement triviale pour $X_t$ d\'ecrivant la famille
d'hypersurfaces param\'etr\'ee par $D_i$.

Si $\delta_i\neq 0$, le morphisme 
$\beta_{B_i,X_i}:\HH_{N}(B_i^{X_i})\to\HH_{N}(X_i)$ est injectif, et
comme $\bord_{B_i}:\HH_{N+1}(B_i,B_i^{X_i})\to\HH_{N}(B_i^{X_i})$ est
un isomorphisme, le morphisme compos\'e $\beta_{B_i,X_i}\circ\bord_{i_B}:
\HH_{N+1}(B_i,B_i^{X_i})\to\HH_{N}(X_i)$ est injectif. Nous
avons $\beta_{B_i,X_i}\circ\bord_{B_i}=\bord\circ\beta_{B_i,Y,X_i}$~;
le morphisme transpos\'e est donc le morphisme surjectif
$\gamma_{B_i,Y\setminus X_i}\circ\tube:\HH_{N}(X_i)\to
\HH_{N+1}(\overline{B}_i^{Y\setminus X_i},\Sigma_i^{Y\setminus X_i})$.
En particulier, le noyau de 
$\gamma_{B_i,Y\setminus X_i}:\HH_{N+1}(Y\setminus X_i)\to
 \HH_{N+1}(\overline{B}_i^{Y\setminus X_i},\Sigma_i^{Y\setminus X_i})$ et
l'image de $\tube:\HH_{N}(X_i)\to\HH_{N+1}(Y\setminus X_i)$ engendrent 
$\HH_{N+1}(Y\setminus X_i)$~: pour tout
$\Lambda\in\HH_{N+1}(Y\setminus X_i)$ il existe 
$\Lambda_{\gamma}\in\kker\,\gamma_{B_i,Y\setminus X_i}$ et
$\lambda\in\HH_{N}(X_i)$ tels que 
$\Lambda=\Lambda_{\gamma}+\tube\,\lambda$. Comme la classe
$\Lambda_{\gamma}$ est \`a support dans $\Gamma_i^{Y\setminus X_i}$,
elle n'intersecte pas $\Delta_i$, et comme la famille 
$(\Gamma_i^{Y\setminus X_t})$ est topologiquement triviale pour $X_t$
d\'ecrivant la famille d'hypersurfaces param\'etr\'ee par $D_i$, la
classe $\Lambda_{\gamma}$ est $g_i$-invariante~; la formule du lemme
est donc vraie pour $\Lambda_{\gamma}$. D'autre part, comme
$\bord\,\Delta_i=\delta_i$ et que bord est la transpos\'ee de tube, la
formule du lemme pour $\tube\,\lambda$ r\'esulte imm\'ediatement de la
formule de Picard-Lefschetz usuelle. 
\cucu

\medskip

\noindent{\it Remarque.\ ---\ }
Il existe une formule duale pour l'homologie relative~: pour tout
$\Lambda\in\HH_{N+1}(Y,X_i)$ et nous avons $g_i(\Lambda)=\Lambda+
\varepsilon\langle\lambda\,|\,{\rm tube}\,\delta_i\rangle\,\Delta_i$
avec $\varepsilon=1$ ou $-1$.

\medskip

\noindent
{\it Fin de la preuve de la 
proposition~\ref{monodromie-sur-complementaire}.\ ---\ }
Soit $\Lambda\in\HH_{N+1}(Y\setminus X)\setminus\{0\}$. 

Pour tout 
$i\in\{1,\dots, r\}$ nous notons abusivement $\Delta_i\in\HH_{N}(Y,X)$
la pr\'eimage de $\Delta_i\in\HH_{N}(Y,X_i)$ par transport plat le long du 
chemin $L_i$. D'apr\`es \cite{cours}, proposition 2.27, les cycles
$\Delta_i$, $i\in\{1,\dots, r\}$ engendrent $\HH_{N+1}(Y,X)$. Comme la
forme d'intersection induit une dualit\'e parfaite
$\HH_{N+1}(Y,X)\simeq\HH_{N+1}(Y\setminus X)^{\vee}$, il  
existe donc un entier $i\in\{1,\dots, r\}$ tel que 
$\langle\Lambda\,|\,\Delta_i\rangle\neq 0$.
D'apr\`es la formule de Picard-Lefschetz pour le compl\'ementaire,
nous avons alors $\tube\,\delta_i=
\frac{\varepsilon}{\langle\Lambda|\Delta_i\rangle}
(g_i(\Lambda)-\Lambda)$~; la classe $\tube\,\delta_i$
appartient donc \`a la repr\'esentation engendr\'ee par $\Lambda$. 

Comme l'action de $\pi_1(D^*,F)$ sur
$\HH_{N}(X)_{\ev}$ est irr\'eductible et commute avec l'application
tube, la $\pi_1(D^*,F)$-repr\'esentation 
engendr\'ee par $\Lambda$ contient $\tube\,\HH_{N}(X)^{\ev}$. La
proposition~\ref{monodromie-sur-complementaire} r\'esulte alors de ce
que l'espace $\C\Lambda+\tube\,\HH_{N}(X)^{\ev}$ 
est $\pi_1(D^*,F)$-invariant.
\cucu

\subsubsection{\'Etude de 
$\HH_N(X\setminus(X_K\setminus C_K))$}\label{referee}

Dans cette section nous montrons un r\'esultat pr\'eliminaire qui d\'ecrit
l'homologie d'une vari\'et\'e intervenant dans la 
proposition~\ref{decompositionpremiere}. Il sert dans les
sections~\ref{monosurface} et~\ref{genrec}, o\`u nous nous 
trouverons dans l'une des situations d\'ecrites ci-dessous.

\medskip

\noindent{\it Situation g\'en\'erale (sections~\ref{genrec}
et~\ref{preuveaux}) \ ---\ } 
Soient $N,\ d$ et $e$ des entiers strictement positifs tels que 
$N\geq 2$ et $e<d$, et soit $Y$ une vari\'et\'e projective lisse de
dimension $N+1$ munie d'un faisceau inversible tr\`es ample $\anneau
(1)$. Pour tout $i\in\NN^*$ et $S\in\HH^0(Y,\anneau(i))\setminus\{0\}$
nous notons $X_S$ l'hypersurface associ\'ee et $j_S:X_S\to Y$
l'immersion ferm\'ee. 
Si $X_A$, $X_K$ et $X_Q$ sont des hypersurfaces de $Y$, nous posons
$X_{A,K}=X_A\cap X_K$ et $X_{A,K,Q}=X_A\cap X_K\cap X_Q$~; nous notons
$j^{A}_{A,K}:X_{A,K}\to X_A$, $j^{A,K}_{A,K,Q}:X_{A,K,Q}\to X_{A,K}$, 
$j^{A}_{A,K,Q}:X_{A,K,Q}\to X_A$, $j_{A,K}:X_{A,K}\to Y$ 
et $j_{A,K,Q}:X_{A,K,Q}\to Y$ les immersions ferm\'ees.

Fixons $A\in\HH^0(Y,\anneau(i))\setminus\{0\}$ tel que $X_A$ est 
lisse.  Notons ${\cal V}^i\subset\HH^0(Y,\anneau(i))$,
et ${}^A{\cal V}^i\subset\HH^0(X_A,\anneau(i))$ les ouverts param\'etrant 
(\`a la multiplication par un scalaire pr\`es) les hypersurfaces lisses 
de $Y$ et de $X_A$ de classe $c_1(\anneau (i))$. 
Fixons ensuite $Q\in\HH^0(Y,\anneau(d))$ en intersection 
compl\`ete avec $A$. Notons ${\cal V}^{i}_{A,Q}\subset\HH^0(Y,\anneau(i))$
l'espace des \'el\'ements $K$ en intersection compl\`ete avec $A$ et $Q$ et
tels que les vari\'et\'es $X_K$, $X_{A,K}$ et $X_{A,K,Q}$
sont lisses. Notons ${}^A{\cal V}_Q^{i}\subset\HH^0(X_A,\anneau(i))$
l'espace des \'el\'ements $K$ en intersection 
compl\`ete avec $Q_{|X_A}$ tels que les vari\'et\'es 
$X_{A,K}$ et $X_{A,K,Q}$ sont lisses. 
Remarquons que les espaces ${\cal V}^{i}_{A,Q}$ et 
${}^A{\cal V}^{i}_{Q}$ sont non vides si et seulement si la vari\'et\'e 
$X_{A,Q}$ a au plus des singularit\'es isol\'ees. 
Le morphisme de restriction 
${}^A\psi^i:\HH^0(Y,\anneau(i))\to\HH^0(X_A,\anneau(i))$ envoie 
${\cal V}^{i}_{A,Q}$ sur ${}^A{\cal V}_Q^i\cap\im{}^A\psi^i$.

\medskip

\noindent{\it Situation particuli\`ere (section~\ref{monosurface}).\ ---\ }
Soit $d$ un entier strictement positif, soit $S$ une surface
projective lisse munie d'un faisceau inversible tr\`es ample 
$\anneau(1)$ et soit $C$ une courbe r\'eduite de $S$. Pour tout
$K\in\HH^0(S,\anneau(d))\setminus\{0\}$ notons $S_K$ la courbe
associ\'ee, posons $C_K=S_K\cap C$ et notons $j^C_{K}:C_K\to C$
l'immersion ferm\'ee. Notons ${\cal V}\subset\HH^0(S,\anneau (d))$
l'ouvert param\'etrant (\`a la 
multiplication par un scalaire pr\`es) les courbes lisses de $S$ de
classe $c_1(\anneau (d))$ et ${\cal V}_C\subset{\cal V}$
l'ouvert des \'el\'ements $K$ tels que le sch\'ema
$C_K$ est lisse de dimension $0$.  
 
Remarquons que si nous posons $S=X_A$, la situation g\'en\'erale dans le
cas $N=2$ permet de retrouver la situation particuli\`ere avec
l'hypoth\`ese suppl\'ementaire que la courbe $C$ est de classe
$c_1(\anneau (d))$.  

\medskip

Nous nous pla\c cons soit dans la situation particuli\`ere et fixons
$K\in{\cal V}_C$, soit dans la situation g\'en\'erale et fixons
$K\in{}^A{\cal V}^{d-e}_{Q}$.
Le lemme~\ref{thom} ci-dessous emprunte ses notations \`a la situation
particuli\`ere~; il est encore vrai dans la situation g\'en\'erale avec les 
notations $S=X_A$, $C=X_{A,Q}$, $S_K=X_{A,K}$, $C_K=X_{A,K,Q}$,
${\cal V}={}^A{\cal V}^{d-e}$ et ${\cal V}_C={}^A{\cal V}^{d-e}_{Q}$. 

\begin{figure}
\centering{
\includegraphics[height=4.5cm]{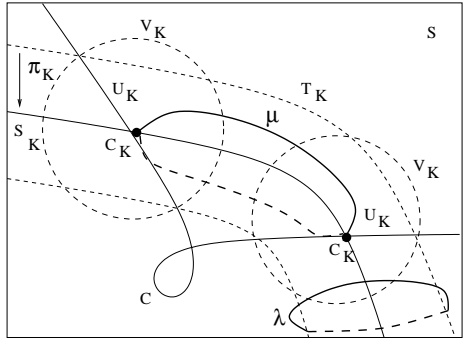}
\caption{Exemples de cycles de 
         $\HH_N(S\setminus(S_K\setminus C_K))$~: 
         $\lambda\in\alpha(\HH_N(S\setminus S_K))$, $\mu\in
         u(\HH_{N-1}(S_K,C_K))$}
\label{cycles}}
\end{figure}

\begin{lemme}\label{thom}
Pour tout $i\in\NN$ nous avons le diagramme suivant de 
$\pi_1({\cal V}_C,K)$-modules. Il est commutatif et ses lignes 
sont des suites exactes. 
$$\xymatrix@C-2mm{
& & \HH_i(C)\dto_{\phi}\drto^{(j^C_{K})^*} & & \\
0\rto&\HH_i(S\setminus S_K)\dto^{\ouv}\rto^{\alpha\qquad}
     &\HH_i(S\setminus(S_K\setminus C_K))\ar@{=}[d]\rto^{\qquad\psi}
     &\HH_{i-2}(C_K)\rto&0\\
     &\HH_i(S)
     &\HH_i(S\setminus(S_K\setminus C_K))\lto_{v\qquad}
     &\HH_{i-1}(S_K,C_K)\lto_{\qquad u}\uto_{\bord}&\HH_{i+1}(S).\lto\\ }$$
Les fl\`eches $\alpha$, $v$ et $\phi$ sont induites par les
immersions. Les fl\`eches $\psi$ et $u$ sont d\'efinies dans la preuve.
Si $\hh_i(S_K)_{\ev}\neq 0$ alors $\im(\phi)\cap\im(\alpha)=0$.
\end{lemme}
Le dessin~\ref{cycles} donne des exemples de cycles de
$\HH_N(S\setminus(S_K\setminus C_K))$.

\medskip

\preuve
Les lignes du diagramme se d\'eduisent des suites exactes longues
d'homologie relative des couples 
\begin{itemize}
\item $(S\setminus(S_K\setminus C_K),S\setminus S_K)$ pour la
      premi\`ere ligne et
\item $(S,S\setminus(S_K\setminus C_K))$ pour la seconde ligne.
\end{itemize}
Nous devons montrer les assertions suivantes.
\begin{enumerate}
\item Nous avons des isomorphismes de $\pi_1({\cal V}_C,K)$-modules
      $\HH_i(S\setminus(S_K\setminus C_K),S\setminus S_K)
      \simeq \HH_{i-2}(C_K)$ et le triangle du diagramme est
      commutatif. \label{thom1} 
\item Nous avons des isomorphismes de $\pi_1({\cal V}_C,K)$-modules
      $\HH_{i+1}(S,S\setminus(S_K\setminus C_K))\simeq
      \HH_{i-1}(S_K,C_K)$ et le rectangle de droite du diagramme est
      commutatif (la commutativit\'e du rectangle de gauche est
      \'evidente).\label{thom2}
\item La premi\`ere ligne est exacte \`a gauche et \`a droite,
      autrement dit les morphismes de liaison
      $\HH_{i-2}(C_K)\to\HH_{i-1}(S\setminus S_K)$ sont nuls~; nous
      distinguons les cas
      \begin{enumerate}
           \item $i\leq\dim S$ et \label{exa1}
           \item $i>\dim S$.\label{exa2}
      \end{enumerate}
\item Si $\hh_i(S_K)_{\ev}\neq 0$ alors
      $\im(\phi)\cap\im(\alpha)=0$. \label{intvide} 
\end{enumerate}

Comme les vari\'et\'es $S_K$ et $C_K$ sont lisses et que la
vari\'et\'e $C$ est lisse au voisinage de $C_K$, il existe
des voisinages tubulaires ouverts $V_K$ de $C_K$ dans $S$, $U_K$ de
$C_K$ dans $C$ et $T_K$ de $S_K$ dans $S$ tels que $U_K=C\cap
V_K=C\cap T_K$ ({\it cf.} dessin~\ref{cycles}). 
La construction de ces voisinages tubulaires pour $C$ et $S_K$
s'\'etend \`a une famille param\'etr\'ee par ${\cal V}_C$, ce qui
assure la compatibilit\'e des constructions ci-dessous avec l'action
de $\pi_1({\cal V}_C,K)$. 

\medskip

\noindent {\it Preuve de l'assertion~\ref{thom1}. \ --- \ } 
Remarquons que $V_K\setminus(S_K\cap V_K)$ se
r\'etracte par d\'eformation sur $U_K\setminus C_K$.
Par  excision, r\'etraction et isomorphisme de
Thom nous avons alors les isomorphismes successifs
\begin{eqnarray*}
  \HH_i(S\setminus(S_K\setminus C_K),S\setminus S_K)&\simeq&
  \HH_i(V_K\setminus((S_K\setminus C_K)\cap V_K),
             V_K\setminus (S_K\cap V_K))\\
   &\simeq&\HH_{i}(U_K,U_K\setminus C_K)\\
   &\simeq&\HH_{i-2}(C_K).
\end{eqnarray*}
La commutativit\'e du triangle du diagramme s'obtient en prenant l'image
r\'eciproque de ces isomorphismes dans $\HH_i(C,C\setminus C_K)$
et en observant que le morphisme
$(j^C_{K})^*:\HH_{i}(C)\to\HH_{i-2}(C_K)$ se factorise par 
$\HH_i(C,C\setminus C_K)\simeq\HH_{i}(U_K,U_K\setminus C_K)$.

\medskip

\noindent {\it Preuve de l'assertion~\ref{thom2}. \ --- \ }
Comparons la suite exacte longue d'homologie relative du couple 
$(S_K,C_K)$ \`a celle du triplet 
$(S,S\setminus(S_K\setminus C_K),S\setminus S_K)$~: pour tout 
$j\in\NN$ nous avons des isomorphismes de Thom
$\HH_{j-2}(S_K)\simeq\HH_j(S,S\setminus S_K)$, et, par 
l'assertion~\ref{thom1}, 
$\HH_j(S\setminus(S_K\setminus C_K),S\setminus S_K)\simeq\HH_{j-2}(C_K)$.
D'autre part, la projection naturelle $\pi_K:T_K\to S_K$ induit le
morphisme 
$$\HH_{i-1}(S_K,C_K)\stackrel{\pi_K^*}{\longrightarrow}
\HH_{i+1}(T_K,T_K\setminus(S_K\setminus C_K))\simeq
\HH_{i+1}(S,S\setminus(S_K\setminus C_K))$$
o\`u le dernier isomorphisme est l'excision. L'isomorphisme
$\HH_{i+1}(S,S\setminus(S_K\setminus C_K))\simeq\HH_{i-1}(S_K,C_K)$
r\'esulte alors du ``lemme des 5''. 

Le rectangle de droite du diagramme est commutatif par cette construction. 
\cucu

\medskip

\noindent {\it Preuve de l'assertion~\ref{exa1}. \ --- \ }
Par la suite exacte longue d'homologie
relative du couple $(S\setminus(S_K\setminus C_K),S\setminus S_K)$
l'assertion~\ref{exa1} est \'equivalente \`a l'injectivit\'e de
$\alpha:\HH_{i-1}(S\setminus S_K)\to
\HH_{i-1}(S\setminus(S_K\setminus C_K))$ pour $i-1\leq\dim S$~; il
suffit donc de montrer que le morphisme 
${\rm ouv}=v\circ\alpha:\HH_{i-1}(S\setminus S_K)\to\HH_{i-1}(S)$
est injectif,
ce qui r\'esulte de ce que son noyau $\HH_{i-2}(S_K)^{\ev}$ est nul
pour $i\leq\dim S$. 
\cucu

\medskip

\noindent {\it Preuve de l'assertion~\ref{exa2}. \ --- \ }
Dans la situation particuli\`ere nous avons $\dim S=2$ et 
$\dim C_K=0$, donc $\HH_{i-2}(C_K)=0$ pour $i>\dim S$ et
l'assertion~\ref{exa2} est vraie.  
Nous nous pla\c cons donc dans la situation g\'en\'erale.

Par construction, les morphismes de liaison
$\HH_{i-2}(C_K)\to\HH_{i-1}(S\setminus S_K)$ se factorisent par
$$\xymatrix@C-2mm{
\HH_{i-2}(C_K)\ar[d]\ar[r]^{{\rm thom}\qquad}& 
  \HH_{i}(U_K,U_K\setminus C_K)\ar[d]\ar[r]&
  \HH_{i}(S\setminus(S_K\setminus C_K),S\setminus S_K)\ar[d]\\
\HH_{i-2}(S_K)\ar[rd]_{\rm tube}\ar[r]^{{\rm thom}\qquad}& 
  \HH_{i}(T_K,T_K\setminus S_K)\ar[d]^{\rm bord }\ar[r]&
  \HH_{i}(S,S\setminus S_K)\ar[ld]^{\rm bord }\\
&\HH_{i-1}(S\setminus S_K)&\\}$$
o\`u les fl\`eches non nomm\'ees sont induites par les immersions.
L'assertion~\ref{exa2} est \'equivalente \`a ce que l'image 
de $\HH_{i-2}(C_K)\to\HH_{i-2}(S_K)$ est incluse dans le noyau de
$\HH_{i-2}(S_K)\stackrel{\tube}{\longrightarrow}
 \HH_{i-1}(S\setminus S_K)$.

Or d'une part nous avons 
$$\kker(\HH_{i-2}(S_K)\stackrel{\tube}{\longrightarrow}
\HH_{i-1}(S\setminus S_K))=
\kker(\HH_{i-2}(S_K)\to\HH_{i-2}(S_K)^{\ev})\ ;$$
d'autre part, comme $C\subset S$ est tr\`es ample, de classe un multiple
rationnel de $c_1({\cal L})$, pour $i>\dim S$
l'application de restriction
$\HH_{i}(S_K)\to\HH_{i-2}(C_K)$ est surjective et nous avons
$${\rm Coker}(\HH_{i-2}(C_K)\to\HH_{i-2}(S_K))=
  {\rm Coker}(\HH_{i}(S_K)\stackrel{c_1({\cal L})\frown}{\longrightarrow}
  \HH_{i-2}(S_K))=\HH_{i-2}(S_K)^{\prim}.$$
L'assertion~\ref{exa2} r\'esulte alors de ce que $\HH_{i-2}(S_K)^{\ev}$ est
un quotient de $\HH_{i-2}(S_K)^{\prim}$. 
\cucu

\medskip

\noindent {\it Preuve de l'assertion~\ref{intvide}. \ --- \ } 
Elle se d\'eduit de l'\'etude de l'action de monodromie 
$\pi_1({\cal V}_C,K)$ sur la premi\`ere ligne du diagramme. La
vari\'et\'e ${\cal V}_C$ est alors un ouvert 
de Zariski non vide de ${\cal V}$, le morphisme
$\pi_1({\cal V}_C,K)\to\pi_1({\cal V},K)$ est donc surjectif. Comme il
factorise l'action de $\pi_1({\cal V}_C,K)$ sur $\HH_{i}(S\setminus
S_K)$, d'apr\`es la proposition~\ref{monodromie-sur-complementaire}, si
$\hh_i(S_K)_{\ev}\neq 0$, l'espace $\HH_{i}(S\setminus S_K)$ n'a pas
de classes $\pi_1({\cal V}_C,K)$-invariantes. Comme le morphisme
$\alpha$ est injectif, l'espace $\im(\alpha)$ n'a donc pas de classes 
$\pi_1({\cal V}_C,K)$-invariantes. Or l'espace $\im(\phi)$ est 
$\pi_1({\cal V}_C,K)$-invariant, ce qui implique l'assertion~\ref{intvide}.
\cucu

\subsubsection{Cas des surfaces}\label{monosurface} 

Dans cette section nous nous pla\c cons dans la situation
particuli\`ere de la section~\ref{referee}. Nous faisons
l'hypoth\`ese suppl\'ementaire $d\geq 2$. Pour toute partie
$E\subset\HH_2(S\setminus(S_K\setminus C_K))$ 
nous notons $H(E)\subset\HH_2(S\setminus (S_K\setminus C_K))$
la $\pi_1({\cal V}_C,K)$-repr\'esentation engendr\'ee par $E$.

\begin{prop}\label{surfaces} 
Pour tout
$\lambda\in\HH_2(S\setminus(S_K\setminus C_K))\setminus\im(\phi)$
nous avons $\dim H(\lambda)\geq\hh_1(S_K)_{\ev}$. 
\end{prop}

\preuve 
Nous supposons $\hh_i(S_K)_{\ev}\neq 0$ et nous \'etudions le
comportement de la classe $\lambda$ dans le diagramme du
lemme~\ref{thom}. 
  
\medskip

Si $\lambda\in\im(\alpha)+\im(\phi)$, nous \'ecrivons
$\lambda=\alpha(\mu)+\nu$ avec $\mu\in\HH_2(S\setminus S_K)$ et 
$\nu\in \im(\phi)$. Comme $\alpha$ et $\phi$ sont des morphismes
de $\pi_1({\cal V}_C,K)$-modules et que leurs images sont en somme
directe, nous avons $H(\lambda)=H(\alpha(\mu))\oplus H(\nu)$. 
Comme $\lambda\not\in\im(\phi)$, nous avons 
$\mu\neq 0$. D'apr\`es la
proposition~\ref{monodromie-sur-complementaire}, la classe
$\mu$  engendre une $\pi_1({\cal V}_C,K)$-repr\'esentation de
dimension sup\'erieure ou \'egale \`a $\hh_1(S_K)^{\ev}$
(comme dans la preuve de l'assertion ~\ref{intvide} de la
section~\ref{referee})~; comme $\alpha$ est injectif, nous avons 
$\dim\,H(\alpha(\mu))\geq\hh_1(S_K)_{\ev}$, donc
$\dim\,H(\lambda)\geq\hh_1(S_K)_{\ev}$.

\medskip

Nous supposons $\lambda\not\in\im(\alpha)+\im(\phi)$ et
posons $\psi(\lambda)=\sum_{z\in C_K}\eta_z[z]$, avec $\eta_z\in\C$.
Comme le triangle du diagramme du lemme~\ref{thom} est commutatif, il
existe des points $x\in C_K$ et $y\in C_K$ appartenant \`a la m\^eme composante 
irr\'eductible $C'\subset C$ et tels que $\eta_x\neq\eta_y$.
Comme $\pi_1({\cal V}_C,K)$ se surjecte sur le produit des groupes de
permutations des points de $C_K$ appartenant \`a une m\^eme composante
irr\'eductible de $C$ ({\it cf}.~\cite{acgh}, p. 111 pour le cas o\`u $C$
est irr\'eductible~; le cas g\'en\'eral se d\'emontre de la m\^eme mani\`ere), il
existe un \'el\'ement 
$g\in\pi_1({\cal V}_C,K)$ dont l'action sur $\HH_0(C_K)$ permute $x$ et
$y$ et laisse fixe les autres points de $C_K$. Nous avons donc
$$\psi(g(\lambda)-\lambda)=
g(\psi(\lambda))-\psi(\lambda)=
(\eta_y-\eta_x)([x]-[y]).$$
Posons $\lambda'=\frac{g(\lambda)-\lambda}{\eta_y-\eta_x}$~: nous avons
$g(\lambda')=[x]-[y]$. Comme $\lambda'\in H(\lambda)$,
il suffit de montrer la proposition~\ref{surfaces} pour $\lambda'$.

Comme l'action de $\pi_1({\cal V}_C,K)$ sur $\HH_2(S)$
est triviale, nous avons
$$v(\lambda')=v\left(\frac{g(\lambda)-\lambda}{\eta_y-\eta_x}\right)
=\frac{g(v(\lambda))-v(\lambda)}{\eta_y-\eta_x}=0.$$
D'apr\`es la seconde ligne du diagramme du lemme~\ref{thom},
il existe donc une classe $\mu\in\HH_{1}(S_K,C_K)$ telle que
$u(\mu)=\lambda'$. Comme l'action de $\pi_1({\cal V}_C,K)$ sur
l'image de $\HH_3(S)$ dans $\HH_{1}(S_K,C_K)$ est triviale,
son intersection avec l'image de $\HH_1(S_K)_{\ev}$ dans
$\HH_{1}(S_K,C_K)$ est nulle~; la
restriction de $u$ \`a l'image de $\HH_1(S_K)_{\ev}$ est donc
injective. Il suffit donc de montrer que la
$\pi_1({\cal V}_C,K)$-repr\'esentation engendr\'ee par 
$\mu$ dans $\HH_{1}(S_K,C_K)$ contient l'image de $\HH_1(S_K)_{\ev}$.

Comme le rectangle de droite du diagramme du lemme~\ref{thom} est
commutatif, le bord de  $\mu$ dans $\HH_0(C_K)$ s'identifie \`a 
$\psi(u(\mu))=[x]-[y]$~: c'est un \'el\'ement
de $\HH_0(C'_K)\subset\HH_0(C_K)$, o\`u nous avons pos\'e
$C'_K=C'\cap S_K$. La classe $\mu$ appartient donc au sous-espace
$\HH_{1}(S_K,C'_K)$ de $\HH_{1}(S_K,C_K)$. La
proposition~\ref{surfaces} r\'esulte donc du lemme suivant, 
appliqu\'e \`a $C=C'$.
 \cucu

\begin{lemme}\label{12}
Supposons la courbe $C$ irr\'eductible. Soient $x$ et $y$ des points
distincts de $C_{K}$ et soit $\mu\in\HH_1(S_K,C_K)$ une classe
ayant pour bord $[x]-[y]\in\HH_0(C_K)$. Alors
la $\pi_1({\cal V}_C,K)$-repr\'esentation engendr\'ee par $\mu$
contient l'image de $\HH_1(S_K)_{\ev}$ dans $\HH_{1}(S_K,C_K)$.
\end{lemme}

\preuve
Notons ${\cal P}$ l'hypersurface param\'etrant les courbes $S_L$, 
$L\in\HH^0(S,\anneau(d))$ singuli\`eres et ${\cal Q}$ l'hypersurface 
param\'etrant les courbes $S_L$ qui ne rencontrent pas $C$
transversalement. Nous avons 
${\cal V}_C=\HH^0(S,\anneau(d))\setminus({\cal P}\cup{\cal Q})$.
La preuve du lemme~\ref{12} repose sur l'\'etude de la monodromie locale
de $\pi_1({\cal V}_C,K)$ sur $\HH_1(S_K,C_K)_{\ev}$ au voisinage d'un
point (g\'en\'erique) de ${\cal P}\cap{\cal Q}$. 

Plus pr\'ecis\'ement, soit $z_O\in C$ un point lisse. Comme $d\geq 2$, il
existe un point $O\in {\cal P}\cap{\cal Q}$ tel que
\begin{itemize}
  \item $z_O$ est un point double ordinaire de $S_O$ et
        $S_O\setminus\{z_O\}$ est lisse~; 
  \item $z_O$ est un point de multiplicit\'e $2$ de $C_O$ et 
        $C_O\setminus\{z_O\}$ est lisse.
\end{itemize}
Nous utilisons les notations de la section~\ref{rappelmonodromie}.
D'apr\`es la th\'eorie de Morse, il existe un ouvert diff\'eomorphe \`a
une boule ouverte $B_O$ de centre $z_O$ dans $S$ et une boule ouverte
${\cal D}_O$ de centre $K_O$ dans $\HH^0(S,\anneau(d))$ tels que
\begin{enumerate}
\item\label{hh1} l'hypersurface ${\cal P}\cap{{\cal D}_O}\subset 
                 {\cal D}_O$ est lisse et 
  \begin{itemize}
    \item pour tout $L\in{\cal P}\cap{{\cal D}_O}$ la vari\'et\'e 
      $B_O^{S_L}$ est diff\'eomorphe \`a la r\'eunion de deux disques qui se
      coupent transversalement en l'unique point singulier $z_L$ de
      $S_L$~; chaque disque contient un point de $C_L$ distinct de $z_L$
      pour $L\neq O$~;  
    \item pour tout $L\in{{\cal D}_O}\setminus({\cal P}\cap{{\cal
      D}_O})$ la vari\'et\'e $B_O^{S_L}$ est diff\'eomorphe \`a un cylindre~;  
  \end{itemize}
\item\label{hh2} l'hypersurface ${\cal Q}\cap{{\cal D}_O}\subset
                 {\cal D}_O$ est lisse, rencontre l'hypersurface
                 ${\cal P}\cap{{\cal D}_O}$ uniquement en $O$ et 
  \begin{itemize}
    \item pour tout $L\in{\cal Q}\cap{{\cal D}_O}$ le sch\'ema $B_O^{C_L}$
       est un point double~;
    \item pour tout $L\in{{\cal D}_O}\setminus({\cal Q}\cap{{\cal
       D}_O})$ le sch\'ema $B_O^{C_L}$ est la r\'eunion de deux points lisses~; 
  \end{itemize}
\item\label{hh3} les vari\'et\'es $\bigcup_{L\in {\cal D}_O}B_O^{S_L}$ et 
                 $\bigcup_{L\in {\cal D}_O}B_O^{C_L}$ sont lisses~;
\item\label{hh4} les familles de vari\'et\'es \`a bord
                 $(\Gamma_O^{S_L},\Sigma_O^{S_L})_{L\in{\cal D}_O}$ et
                 $(\Gamma_O^{C_L},\Sigma_O^{C_L})_{L\in {\cal D}_O}$ 
                 sont topologiquement triviales.
\end{enumerate}

\begin{figure}
\centering{
\parbox[c]{3.7cm}{\centering{
\includegraphics[height=3.5cm]{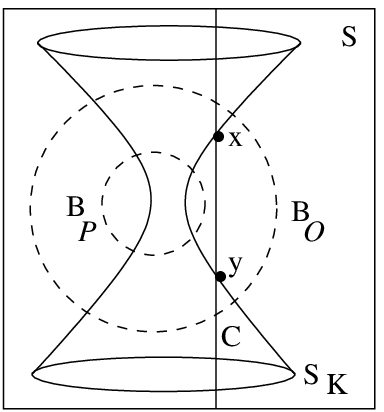} }}
\parbox[c]{3.7cm}{\centering{
\includegraphics[height=3.5cm]{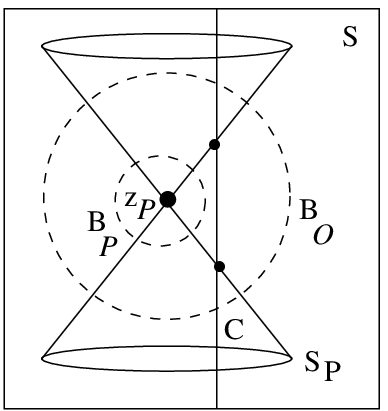} }}
\parbox[c]{3.7cm}{\centering{
\includegraphics[height=3.5cm]{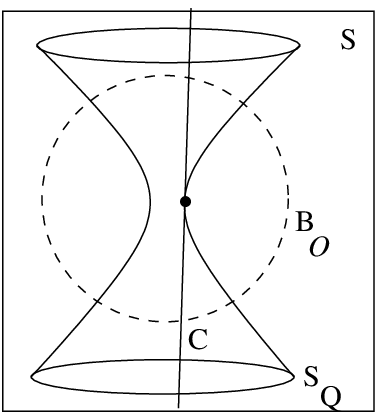} }}
\parbox[c]{3.7cm}{\centering{
\includegraphics[height=3.5cm]{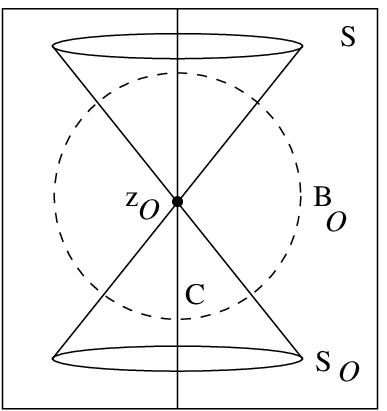} }}
\caption{$K\in{\cal D}_O\setminus(({\cal P}\cup{\cal Q})\cap{\cal D}_O)$,
$P\in{\cal P}\cap{\cal D}_O$, $Q\in{\cal Q}\cap{\cal D}_O$, 
$O={\cal P}\cap{\cal Q}\cap{\cal D}_O$}
\label{X_K} }
\end{figure}

Le dessin~\ref{X_K} repr\'esente les quatre
types de surfaces param\'etr\'ees par ${\cal D}_O$.

\medskip

Comme ${\cal V}_C$ est connexe et que $\pi_1({\cal V}_C,K)$ se
surjecte sur le groupe de permutations des points de $C_K$, il suffit
de montrer le lemme~\ref{12} pour n'importe quel $K\in{\cal V}_C$ et
n'importe quel couple de points dictincts de $C_K$~: nous supposons
$K\in D_O$ et $\{x,y\}=B_O^{C_K}$. Le lemme~\ref{12} r\'esulte alors
des sous-lemmes~\ref{12P} et~\ref{12Q} ci-dessous. 
\cucu

\begin{sous-lemme}\label{12P}
Le lemme~\ref{12} est vrai pour une classe $\nu\in\HH_1(S_K,C_K)$ 
\`a support dans $B_O^{S_K}$ dont le bord est $[x]-[y]$.
\end{sous-lemme}

\begin{figure}
\centering{
\parbox[c]{5.6cm}{\centering{
\includegraphics[height=4.5cm]{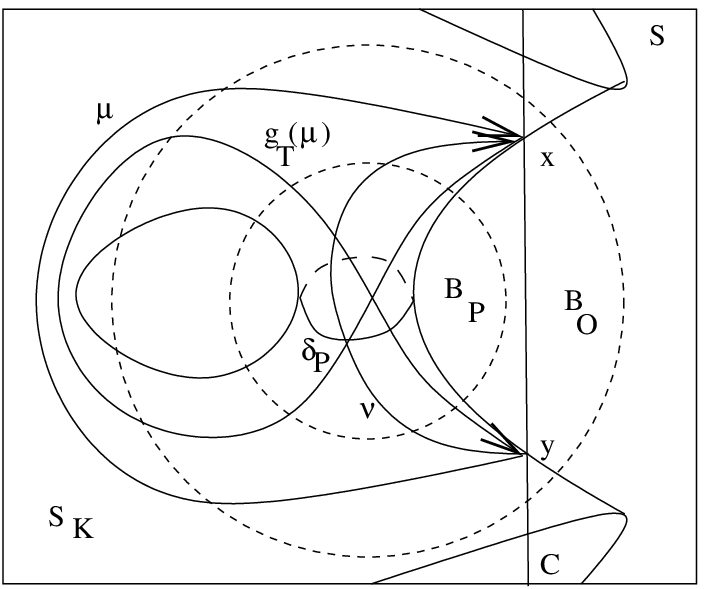}}
\caption{Monodromies}\label{dessinmono} }
\parbox[c]{4.8cm}{\centering{
\includegraphics[height=4.5cm]{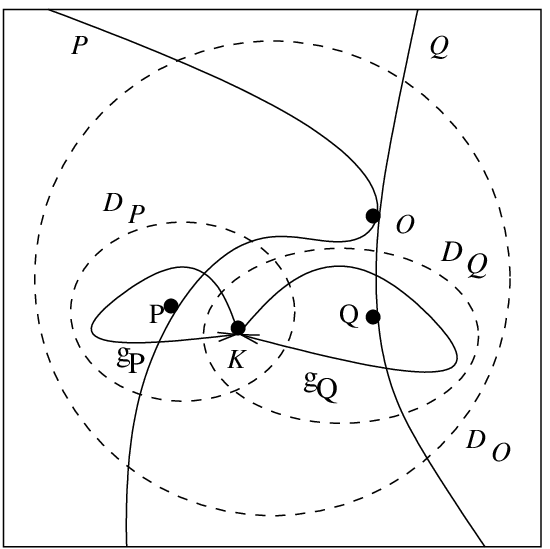}}
\caption{${\cal D}_O$ vu de haut}\label{D_O} }
\parbox[c]{4.7cm}{\centering{
\includegraphics[height=4.5cm]{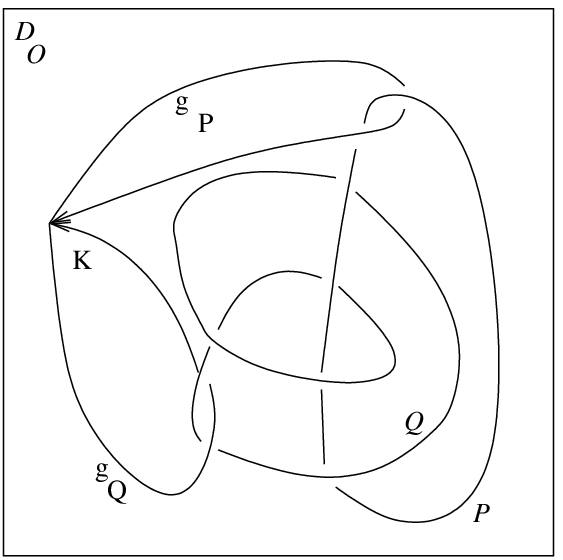}}
\caption{${\cal D}_O$ vu de $O$}\label{pi1} } }
\end{figure}

\preuve
Fixons $P\in {\cal P}\cap{{\cal D}_O}\setminus\{O\}$. 
Comme le point singulier $z_P$ de $S_P$ n'appartient pas \`a $C$, 
il existe un ouvert diff\'eomorphe \`a une boule ouverte $B_P\subset
B_O$ de centre $z_P$ qui ne rencontre pas $C$, tel que la vari\'et\'e
$B_P^{S_P}$ est diff\'eomorphe \`a la r\'eunion de deux disques se
coupant transversalement en $z_P$. D'apr\`es la th\'eorie de Morse, il
existe un disque ${\cal D}_P$ de centre $P$ dans  
${{\cal D}_O}\setminus({\cal Q}\cap{{\cal D}_O})$ qui ne rencontre
${\cal P}$ qu'en $P$, et tel que la vari\'et\'e $\bigcup_{L\in {\cal
D}_P}B_P^{S_L}$ est lisse et la famille de vari\'et\'es \`a bord
$(\Gamma_P^{S_L},\Sigma_P^{S_L})_{L\in {\cal D}_P}$ est
topologiquement triviale ({\it cf.} dessin~\ref{D_O}). 

Sans perdre en g\'en\'eralit\'e, nous supposons 
$K\in{{\cal D}_P}\setminus\{P\}$, en prenant pour $x$ et $y$ les
deux points de $C_K\cap B_O$.  
Comme les familles $(B_O^{S_L}\setminus B_P^{S_L})_{L\in {\cal D}_P}$
et $(B_O^{C_L}\setminus B_P^{C_L})_{L\in {\cal D}_P}$ sont 
topologiquement triviales, l'espace $B_O^{S_K}\setminus B_P^{S_K}$ est
hom\'eomorphe \`a $B_O^{S_O}\setminus B_P^{S_O}$ qui s'identifie \`a
la r\'eunion disjointe de deux cylindres contenant chacun un point de
$C_O$. Donc $B_O^{S_K}\setminus B_P^{S_K}$ est hom\'eomorphe \`a 
la r\'eunion disjointe de deux sous-cylindres du cylindre $B_O^{S_K}$
contenant l'un le point $x$ et l'autre le point $y$ et s\'epar\'es par le
sous-cylindre $B_P^{S_K}$ ({\it cf.} dessins~\ref{X_K}
et~\ref{dessinmono}). L'image $\nu^{B_P,\Sigma_P}$ de $\nu$ dans
$\HH_1(\overline{B}_P^{S_K},\Sigma_P^{S_K})$ est alors la classe d'une
g\'en\'eratrice du cylindre $B_P^{S_K}$. 

Soit $\delta_P^{B_P}\in\HH_1(B_P^{S_K})$ la classe d'un cercle
sur le cylindre $B_P^{S_K}$ qui engendre son groupe fondamental, dont
l'image dans $\HH_1(S_K)$ est le cycle \'evanescent $\delta_P$.
Nous avons $\langle\nu^{B_P,\Sigma_P},\delta_P^{B_P}\rangle=\varepsilon$,
avec $\varepsilon=1$ ou $-1$ selon l'orientation de $\delta_P^{B_P}$
({\it cf.} dessin~\ref{dessinmono}). 

La formule de Picard-Lefschetz pour $\HH_1(S_K)$ s'\'etend \`a
$\HH_1(S_K,C_K)$~: pour tout $\lambda\in\HH_1(S_K,C_K)$, si
$\lambda^{B_P,\Sigma_P}$ d\'esigne l'image de $\lambda$ dans  
$\HH_1(\overline{B}_P^{S_K},\Sigma_P^{S_K})$, si $\Delta_P$ d\'esigne 
l'image de $\delta_P$ dans $\HH_1(S_K,C_K)$, 
et si $g_{S}\in\pi_1({\cal V}_C,K)$ est la classe d'un lacet d'origine
$K$ faisant le tour de $P$ dans ${\cal D}_P$, bien d\'efinie au signe pr\`es,  
nous avons  
$$g_P(\lambda)-\lambda=
  \varepsilon\langle\lambda^{B_P,\Sigma_P},\delta_P^{B_P}\rangle\Delta_P, 
  \quad {\rm avec}\quad \varepsilon=1\ {\rm ou}\ -1.$$
({\it Preuve~:} comme $\Gamma_P^{S_K}$ est connexe par arcs et contient
$C_K$, le cycle $\lambda$ est la somme d'un cycle dans l'image de 
$\HH_1(\Gamma_P^{S_K},C_K)$ et d'un cycle dans $\rel(\HH_1(S_K))$. 
Le premier cycle est $g_P$-invariant et pour le second la formule ci-dessus
r\'esulte directement de la formule de Picard-Lefschetz usuelle.)

En particulier, nous avons $g_P(\nu)-\nu=\varepsilon\Delta_P$ avec 
$\varepsilon=1$ ou $-1$. Mais $\Delta_P$ appartient \`a l'image de
$\HH_1(S_K)_{\ev}\setminus\{0\}$, et comme l'action de 
$\pi_1({\cal V}_C,K)$ sur $\HH_1(S_K)_{\ev}$ est irr\'eductible,
la $\pi_1({\cal V}_C,K)$-repr\'esentation engendr\'ee par $\nu$ contient
l'image de $\HH_1(S_K)_{\ev}$. 
\cucu

\begin{sous-lemme}\label{12Q}
Pour tout $\mu\in\HH_1(S_K,C_K)$, il existe une classe $\nu\in H(\mu)$
\`a support dans $B_O^{S_K}$  dont le bord est $[x]-[y]$. 
\end{sous-lemme}

\preuve
Il existe un espace diff\'eomorphe \`a un disque 
${\cal D}_Q\subset{{\cal D}_O}\setminus({\cal P}\cap {\cal D}_O)$ qui
contient $K$, rencontre ${\cal Q}$ en un unique point $Q$, et tel que 
la vari\'et\'e $\bigcup_{L\in {\cal D}_Q}B_O^{C_L}$ est lisse. Soit
$g_{Q}\in\pi_1({\cal V}_C,K)$ la classe d'un lacet $\Gamma_Q$
d'origine $K$ faisant le tour de $Q$ dans ${\cal D}_Q$, bien d\'efinie
au signe pr\`es. 

Nous \'etudions d'abord l'action de monodromie de $g_Q$ sur
$\HH_0(C_K)$. Comme la vari\'et\'e $\bigcup_{L\in {\cal
D}_Q}B_O^{C_L}$ est un rev\^etement double {\em lisse} de ${\cal D}_Q$,
ramifi\'e en $Q$, la vari\'et\'e $\bigcup_{L\in\Gamma_Q}B_O^{C_L}$ est
un rev\^etement double {\em connexe} de $\Gamma_Q$. Nous avons donc
$g_Q([x]-[y])=[y]-[x]$ (nous pouvons aussi d\'eduire cette \'egalit\'e de la
formule de Picard-Lefschetz appliqu\'ee \`a $[x]-[y]$, en remarquant que
$[x]-[y]$ est le cycle \'evanescent).  

Nous \'etudions ensuite l'action de $g_Q$ sur $\HH_1(S_K,C_K)$. Nous
posons $\nu=\frac{\mu-g_Q(\mu)}{2}$ ({\it cf.}
dessin~\ref{dessinmono}). Alors 
$$\bord(\nu)=\bord\left(\frac{\mu-g_Q(\mu)}{2}\right)
            =\frac{\mu-g_Q(\bord(\mu))}{2}
            =\frac{[x]-[y]-g_Q([x]-[y])}{2}=[x]-[y]. $$
D'autre part, l'action de $g_Q$ laisse invariante l'image de $\mu$ dans
 $\HH_1(S_K,B_O^{S_K}\cup C_K)$~; le cycle 
$\nu=\frac{\mu-g_Q(\mu)}{2}$ appartient donc \`a l'image de 
$\HH_1(B_O^{S_K},\{x,y\})$ dans $\HH_1(S_K,C_K)$. 
Comme $\nu\in H(\mu)$, ceci ach\`eve la preuve du sous-lemme~\ref{12Q}.
\cucu 

\medskip

\noindent{\it Remarque sur le lemme~\ref{12}.\ ---\ }
Les hypersurfaces ${\cal P}\cap {\cal D}_O$ et ${\cal Q}\cap {\cal
D}_O$ ne se coupent pas transversalement, un calcul direct montre
en effet qu'elles sont tangentes. Ceci explique que les actions de
$g_P$ et $g_Q$ ne commutent pas. Plus pr\'ecis\'ement, nous pouvons
montrer que pour un choix convenable de $g_P$ et $g_Q$,
le groupe $\pi_1({\cal D_O},K)$ est le groupe libre engendr\'e par $g_P$
et $g_Q$ quotient\'e par la relation $(g_Pg_Q)^2=(g_Qg_P)^2$ 
(le dessin~\ref{pi1} repr\'esente alors l'intersection de ${\cal D}_O$
avec une petite sph\`ere de centre $O$). Sa repr\'esentation de monodromie 
sur $\HH_1(S_K,C_K)$ passe au quotient par la relation
suppl\'ementaire $g_Q^2=1$ et la repr\'esentantion quotient 
est fid\`ele. Ceci m\`ene \`a une autre preuve du lemme~\ref{12},
qui n'est pas plus simple.

\medskip

\noindent{\it Remarque sur la proposition~\ref{surfaces}.\ ---\ }
La proposition~\ref{surfaces} reste vraie dans la situation
g\'en\'erale de la section~\ref{referee}. La preuve est identique \`a
celle de la situation particuli\`ere (il faut remplacer le cycle
$[x]-[y]$ par une sph\`ere \'evanescente), sauf pour  la
surjectivit\'e de  $\pi_1({\cal V}_C,K)$ sur le produit des groupes de
permutations des points de $C_K$ appartenant \`a une m\^eme composante
irr\'eductible de $C$ qui doit \^etre remplac\'ee pour $N>2$ par
l'irr\'eductibilit\'e de l'action de monodromie de 
$\pi_1({}^A{\cal V}^{d,d-e}_{Q},K)$ sur  
$\HH_{N-2}(X_{A,K,Q}/\im(j_{A,K,Q}^{A,Q})^*)$, qui est fausse pour
$N=2$ et p\'enible \`a d\'emontrer pour $N>2$. \'Enoncer la
proposition~\ref{surfaces} pour la situation g\'en\'erale ne simplifie
pas la preuve du th\'eor\`eme~\ref{mono}~; nous avons donc choisi
de l'\'enoncer pour la situation particuli\`ere qui privil\'egie
l'intuition g\'eom\'etrique.  

\subsection{Monodromie locale d'une famille d'hypersurfaces au
voisinage d'une d\'eg\'en\'erescence en la r\'eunion de deux
hypersurfaces}\label{genrec} 

Nous nous pla\c cons dans la situation
g\'en\'erale de le section~\ref{referee}. 
Nous fixons $A\in{\cal V}^{e}$ et $Q\in\PP\HH^0(Y,\anneau(d))$ tels
que
$X_{A,Q}$ a au plus des singularit\'es isol\'ees si bien que 
${\cal V}^{d-e}_{A,Q}$ est non vide.

\subsubsection{La filtration dissym\'etrique}\label{rec}   

Comme pour tout $K\in{\cal V}^{d-e}_{A,Q}$
la sous-vari\'et\'e de $Y\times\PP_{\C}^1$ d\'efinie 
par l'\'equation $AK+tQ=0$, $t\in\PP_{\C}^1$ est lisse pour $t\neq 0$,
pour $t\in\PP_{\C}^1$ g\'en\'erique l'hypersurface $X_{AK+tQ}$ est lisse.
Soit $\rho_{K}$ le plus petit r\'eel (\'eventuellement infini) tel que
si $\Delta\subset\C$ est un disque ouvert de rayon $\rho_{K}$ et de
centre~$0$, alors pour tout $t\in\Delta\setminus\{0\}$ l'hypersurface
$X_{AK+tQ}$ est lisse. Nous posons 
$${\cal U}_{A,Q}=\{K\in {\cal V}^{d-e}_{A,Q}\ |\ \rho_{K}>1\}$$
si bien que pour tout $K\in {\cal U}_{A,Q}$ l'hypersurface
$X_{AK+Q}$ est lisse. L'application continue
${\cal U}_{A,Q}\to{\cal V}^{d}$, $K\mapsto AK+Q$ induit un morphisme de
groupes $\pi_1({\cal U}_{A,Q},K)\to\pi_1({\cal V}^{d},AK+Q)$~;  
le groupe $\pi_1({\cal U}_{A,Q},K)$ agit donc par monodromie sur
$\HH_N(X_{AK+Q})$. 

\begin{prop}\label{decompositionpremiere}
Pour tout $K\in {\cal U}_{A,Q}$ nous avons le diagramme suivant, dont
la ligne et les colonnes sont des suites exactes de
$\pi_1({\cal U}_{A,Q},K)$-modules.
$$\xymatrix{
 &0\ar[d]& &0\ar[d]&\\
   &\HH_N(X_A\setminus X_{A,K})\ar[d]^{\alpha}& 
   &\HH_N(X_K)^{\prim}\ar[d]^{\rel}&\\
 0\ar[r]&
 \HH_N(X_A\setminus(X_{A,K}\setminus X_{A,K,Q}))
    \ar[r]^{\qquad F_2}\ar[d]^{\psi}&
 \HH_N(X_{AK+Q})\ar[r]^{F_4/F_2}&
 \HH_N(X_K,X_{A,K})\ar[r]\ar[d]^{\bord}&0\\
   &\HH_{N-2}(X_{A,K,Q})\ar[d]& 
   &\HH_{N-1}(X_{A,K})_{\ev}\ar[d]&\\
 &0& &0&\\}$$
Il existe donc une filtration $F_{\bullet}\HH_N(X_{AK+Q})$ de 
$\pi_1({\cal U}_{A,Q},K)$-modules v\'erifiant 
\begin{enumerate}
\item $F_1\HH_N(X_{AK+Q})=\HH_N(X_A\setminus X_{A,K})$~;
\item $F_2/F_1\HH_N(X_{AK+Q})=\HH_{N-2}(X_{A,K,Q})$~;
\item $F_3/F_2\HH_N(X_{AK+Q})=\HH_N(X_K)^{\prim}$~;
\item $F_4/F_3\HH_N(X_{AK+Q})=\HH_{N-1}(X_{A,K})_{\ev}$~;
\end{enumerate}
et $F_4\HH_N(X_{AK+Q})=\HH_N(X_{AK+Q})$.

De plus, les diagrammes suivants sont commutatifs.
\begin{eqnarray}\xymatrix{
&\HH_N(X_{A,Q})\ar[ld]_{(j^{A,Q}_{A,K,Q})^*}
                   \ar[d]^{\phi}\ar[rd]^{(j^{A,Q}_{AK+Q})_*}&\\
  \HH_{N-2}(X_{A,K,Q})&
  \HH_N(X_A\setminus(X_{A,K}\setminus X_{A,K,Q}))\ar[r]^{\qquad F_2}
                                                 \ar[l]_{\psi\qquad}&
  \HH_N(X_{AK+Q}) }\label{d1}\end{eqnarray}

\begin{eqnarray}\xymatrix{
  \HH_{N+2}(Y)\ar[r]^{j_K^*}\ar[d]^{j_{AK+Q}^*}&\HH_{N}(X_K)\ar[d]^{\rel}\\
  \HH_{N}(X_{AK+Q})\ar[r]^{F_4/F_2}&\HH_{N}(X_K,X_{A,K})}\label{d2}
\end{eqnarray}
\end{prop}

\medskip

\noindent{\it Remarque.\ ---\ }
Nous pouvons d\'ecomposer $\HH_N(X_A\setminus X_{A,K})$ gr\^ace \`a la
suite exacte (\ref{ouvert}). Nous obtenons alors une filtration \`a 5
termes de $\HH_N(X_{AK+Q})$ dont les quotients successifs sont
$$\HH_{N-1}(X_{A,K})^{\ev},\
  \HH_N(X_A)_{\prim},\
  \HH_{N-2}(X_{A,K,Q}),\
  \HH_N(X_K)^{\prim},\
  \HH_{N-1}(X_{A,K})_{\ev}.$$
La forme d'intersection induit un isomorphisme
$\HH_N(X_{AK+Q})\simeq\HH_N(X_{AK+Q})^{\vee}$~; la filtration
duale de la filtration ci-dessus est obtenue en permutant $A$ et
$K$. 

\subsubsection{Preuve de la proposition~\ref{decompositionpremiere}~:
\'etude d'une famille semistable}\label{grrr}

Nous fixons $K\in{\cal U}_{A,Q}$ et nous montrons l'existence,
l'exactitude et la commutativit\'e des 
diagrammes donn\'es par la proposition~\ref{decompositionpremiere}~; la
compatibilit\'e avec l'action de $\pi_1({\cal U}_{A,Q},K)$ s'obtient en
mettant notre construction en familles pour $K$ parcourant 
${\cal U}_{A,Q}$.  

\medskip

Soit ${\cal X}\subset Y\times\Delta$ la vari\'et\'e d\'efinie 
par l'\'equation $AK+tQ=0$, $t\in\Delta$. Pour $N=1$
la vari\'et\'e ${\cal X}$ est lisse. Pour $N\geq 2$ elle est singuli\`ere
exactement en $X_{A,K,Q}\times\{0\}$ (et c'est une singularit\'e
ordinaire). Nous d\'esingularisons la vari\'et\'e ${\cal X}$ en imitant (et
en simplifiant) une construction de~\cite{hodgevariationnel},
section~2a, qui traite le cas $Y=\PP_{\C}^3$. 

Soit $\widetilde{\cal X}$ l'\'eclat\'e 
de ${\cal X}$ en l'id\'eal $(K,Q)$. Remarquons que 
l'id\'eal $(K,Q)$ est localement principal en dehors du lieu 
$X_{A,K,Q}\times\{0\}$~; la projection $p:\widetilde{\cal X}\to
{\cal X}$ induit donc un isomorphisme $\widetilde{\cal X}\setminus
p^{-1}(X_{A,K,Q}\times\{0\})\to{\cal X}\setminus (X_{A,K,Q}\times\{0\})$.
Nous pouvons d\'ecrire $\widetilde{\cal X}$ plus explicitement si nous
munissons le fibr\'e en droites projectives 
${\bf P}={\bf P}(\anneau(d-e)\oplus\anneau(d))$ au-dessus de $Y$ de
coordonn\'ees homog\`enes $(U,V)$~: la vari\'et\'e $\widetilde{\cal X}$ est
alors la sous-vari\'et\'e de ${\bf P}\times\Delta$ d\'efinie par les
\'equations $KV+QU=0$ et $AU-tV=0$. 

Nous notons $\pi:\widetilde{\cal X}\to\Delta$ la projection naturelle et,
pour tout $t\in\Delta$, $\widetilde{\cal X}_t$ la fibre en $t$ de $\pi$.  
Les assertions suivantes sont alors des cons\'equences imm\'ediates de
la description de $\widetilde{\cal X}$ ci-dessus~; elles affirment que
la famille $\pi:\widetilde{\cal X}\to\Delta$ est semi-stable. 
\begin{enumerate}
\item La vari\'et\'e $\widetilde{\cal X}$ est lisse.
\item Pour tout $t\in\Delta\setminus\{0\}$ la vari\'et\'e
$\widetilde{\cal X}_t$ est canoniquement isomorphe \`a $X_{AK+tQ}$, en
particulier elle est lisse.\label{X_t}
\item La vari\'et\'e $\widetilde{\cal X}_0$ est la r\'eunion des
vari\'et\'es $\widetilde{X_A}$, d\'efinie comme l'\'eclat\'e de
$X_A$ le long de $X_{A,K,Q}$ (plong\'ee dans ${\bf P}\times\Delta$ par
les \'equations $t=0$, $A=0$, $KV+QU=0$),
et $X_K$ (plong\'ee dans ${\bf P}\times\Delta$ par les \'equations
$t=0$, $U=0$, $K=0$). Ces vari\'et\'es sont 
lisses, se coupent transversalement et leur intersection est
canoniquement isomorphe \`a $X_{A,K}$.\label{AK}
\end{enumerate}
Nous notons $s:X_{A,K}\to\widetilde{X_A}$ l'immersion ferm\'ee
d\'efinie par l'assertion~\ref{AK}. 

\medskip

La preuve de la proposition~\ref{decompositionpremiere} 
repose sur les lemmes~\ref{semistable},~\ref{contraction},~\ref{sg}
et~\ref{sd} ci-dessous. 

\begin{lemme}\label{semistable}
Pour tout $t\in\Delta\setminus\{0\}$ nous avons une suite exacte
longue 
$$\xymatrix@C-5mm{
  \HH_{i+1}(X_K,X_{A,K})\ar[r]&
  \HH_i(\widetilde{X_A}\setminus s(X_{A,K}))
      \ar[r]^{\qquad\ \widetilde{F_2}}&
  \HH_i(\widetilde{\cal X}_t)\ar[r]^{\widetilde{F_4/F_2}\qquad}&
  \HH_i(X_K,X_{A,K})\ar[r]&
  \HH_{i-1}(\widetilde{X_A}\setminus s(X_{A,K})).}$$
De plus, les diagrammes suivants sont commutatifs~:
$$\xymatrix@C-5mm{
  \HH_{i+2}(\widetilde{Y})\ar[r]\ar[d]&
  \HH_{i}(X_K)\ar[d]^{\rel}& & 
  \HH_i(\widetilde{X_{A,Q}})\ar[r]\ar[rd]&
  \HH_{i}(\widetilde{\cal X}_A\setminus s(X_{A,K}))
       \ar[r]^{\qquad\widetilde{F_2}}\ar[d]&
  \HH_{i}(\widetilde{\cal X}_t)\ar[d]\\
  \HH_{i}(\widetilde{X_t})\ar[r]^{\widetilde{F_4/F_2}\qquad}&
  \HH_{i}(X_K,X_{A,K})& & &
  \HH_{i}(\widetilde{X_A})\ar[r]&\HH_{i}(\widetilde{X}) },$$
o\`u $\widetilde{X_{A,Q}}\subset\widetilde{X_A}$ d\'esigne le
transform\'e strict de $X_{A,Q}$, o\`u $\widetilde{Y}\subset{\bf P}$
d\'esigne l'\'eclat\'e de $Y$ en l'id\'eal $(K,Q)$ et o\`u les
fl\`eches non nomm\'ees sont induites par les immersions.
\end{lemme}

\preuve Comme les vari\'et\'es $\widetilde{\cal X}_t$,
$t\in\Delta\setminus\{0\}$ sont hom\'eomorphes, il suffit de montrer
le lemme~\ref{semistable}  pour un seul $t\in\Delta\setminus\{0\}$.

En mettant en famille la construction de la
section~\ref{rappelmonodromie}, nous obtenons un voisinage tubulaire
ouvert $B_{A,K}$ de $s(X_{A,K})\times\{0\}$, dans ${\bf P}$ et un petit
voisinage $\Delta'$ de $0$ dans $\Delta$
tels que 
\begin{itemize}
\item la famille de vari\'et\'es \`a bord 
$\left(\Gamma_{A,K}^{\widetilde{\cal X}_t},
 \Sigma_{A,K}^{\widetilde{\cal X}_t}\right)_{t\in\Delta'}$ est
topologiquement triviale~;
\item la famille 
$\left(B_{A,K}^{\widetilde{\cal X}_t}\right)_{t\in\Delta'\setminus\{0\}}$
est localement triviale (et ses fibres sont des fibrations en cylindres
au-dessus de $X_{A,K}$)~;
\item la vari\'et\'e $B_{A,K}^{\widetilde{\cal X}_0}$ s'identifie \`a la r\'eunion
d'un voisinage tubulaire $T_A$ de $s(X_{A,K})$ dans $\widetilde{X_A}$
et d'un voisinage tubulaire $T_K$ de $X_{A,K}$ dans $X_K$ qui se
coupent le long de $X_{A,K}$~;
la vari\'et\'e $\Sigma_{A,K}^{\widetilde{\cal X}_0}$ s'identifie \`a la r\'eunion
des bords $\Sigma_A$ et $\Sigma_K$ des voisinages tubulaires $T_A$ et $T_K$,
\end{itemize}

Pour tout $t\in\Delta'$ nous avons donc des hom\'eomorphismes
\begin{eqnarray*}
&\Gamma_{A,K}^{\widetilde{\cal X}_t}&\simeq 
\Gamma_{A,K}^{\widetilde{\cal X}_0}\simeq
(\widetilde{X_A}\setminus T_A)\cup(X_K\setminus T_K)
\simeq(\widetilde{X_A}\setminus s(X_{A,K}))\cup(X_K\setminus X_{A,K})
\quad{\rm et}\\
&\Sigma_{A,K}^{\widetilde{\cal X}_t}&\simeq 
\Sigma_{A,K}^{\widetilde{\cal X}_0}\simeq
\Sigma_A\cup \Sigma_K.
\end{eqnarray*}

Notons $\widetilde{\cal X}_{t,A}$ la composante de 
$\Gamma_{A,K}^{\widetilde{\cal
X}_t}$ hom\'eomorphe \`a $\widetilde{X_A}\setminus s(X_{A,K})$ 
et $\Sigma^{\widetilde{\cal X}_t}_A$ son bord~; de m\^eme, notons
$\widetilde{\cal X}_{t,K}$ la composante de $\Gamma_{A,K}^{\widetilde{\cal
X}_t}$ hom\'eomorphe \`a $X_K\setminus X_{A,K}$ et 
$\Sigma^{\widetilde{\cal X}_t}_K$ son bord. Nous avons alors pour 
tout $t\in\Delta'$
\begin{eqnarray*}
\HH_i(\widetilde{\cal X}_{t,A})&=&
   \HH_i(\widetilde{X_A}\setminus s(X_{A,K}))\qquad {\rm et}\\
\HH_i(\widetilde{\cal X}_t,\widetilde{\cal X}_{t,A})&=&
  \HH_i(\widetilde{\cal X}_{t,K},\Sigma^{\widetilde{\cal X}_t}_K)
\qquad {\rm par\ excision}\\&=& 
  \HH_i(\widetilde{\cal X}_{0,K},\Sigma_K)
\qquad {\rm par\ transport\ plat}\\&=& 
  \HH_i(X_K,T_K)\qquad {\rm par\ excision}\\&=&\HH_i(X_K,X_{A,K})
\qquad {\rm par\ retraction.}
\end{eqnarray*}
Gr\^ace \`a ces isomorphismes, la suite exacte du
lemme~\ref{semistable} s'identifie \`a la suite exacte 
longue d'homologie relative du couple $(\widetilde{\cal
X}_t,\widetilde{\cal X}_{t,A})$.
La commutativit\'e des diagrammes est \'evidente.
\cucu

\begin{lemme}\label{contraction} La projection $p_Y:{\bf P}\to Y$
induit une \'equivalence d'homotopies entre
$\widetilde{X_A}\setminus s(X_{A,K})$ et 
$X_A\setminus(X_{A,K}\setminus X_{A,K,Q})$, et donc un isomorphisme
canonique $\HH_i(\widetilde{X_A}\setminus s(X_{A,K}))\simeq 
\HH_i(X_A\setminus(X_{A,K}\setminus X_{A,K,Q})).$
\end{lemme}

\preuve
Soit  $E_A$  le diviseur exceptionnel de
$p_{Y|\widetilde{X_A}}:\widetilde{X_A}\to X_A$. Alors 
la projection $p_{Y|\widetilde{X_A}}$ d\'efinit par restriction  
\begin{itemize}
\item un isomorphisme
$\widetilde{X_A}\setminus(s(X_{A,K})\cup E_A)\to X_A\setminus
X_{A,K}$, restriction de l'isomorphisme 
$\widetilde{X_A}\setminus E_A\to X_A\setminus X_{A,K,Q}$~;
\item un fibr\'e en droites affines 
$E_A\setminus(s(X_{A,K})\cap E_A)\to X_{A,K,Q}$, inclus dans le fibr\'e 
en droites projectives $E_A\to X_{A,K,Q}$.
\end{itemize}
Notons $T$ un voisinage d'un voisinage tubulaire de $X_{A,K,Q}$
dans $X_A$ et $\widetilde{T}$ l'image inverse de $T$ par 
$p_{Y|\widetilde{X_A}}:\widetilde{X_A}\to X_A$. Alors
$\widetilde{T}$ est un fibr\'e sur $X_{A,K,Q}$ de fibre
$\PP^1_{\C}\times\Delta$ (o\`u $\Delta$ d\'esigne un disque 
ouvert de centre $0$). De plus
\begin{itemize}
\item 
$\widetilde{T}\cap\widetilde{X_A}\setminus s(X_{A,K})$ 
s'identifie \`a un fibr\'e sur $X_{A,K,Q}$ de fibre ${\C}\times\Delta$~; 
\item
$T\cap X_A\setminus(X_{A,K}\setminus X_{A,K,Q})$
s'identifie \`a un fibr\'e sur $X_{A,K,Q}$ de fibre l'image de 
${\C}\times\Delta$ par la contraction qui envoie ${\C}\times\{0\}$ sur
un point.
\end{itemize}
Il existe une \'equivalence d'homotobie entre les deux fibr\'es
ci-dessus induisant l'identit\'e sur leur bord (qui s'identifie \`a un
fibr\'e sur $X_{A,K,Q}$ de fibre $\C\times S^1$). 
\cucu

\begin{lemme}\label{sg}
Le morphisme $\HH_{N+1}(X_{AK+Q})\to\HH_{N+1}(X_K,X_{A,K})$ est surjectif.
\end{lemme}

\preuve
Gr\^ace \`a la commutativit\'e du diagramme~(\ref{d2}), il suffit de
montrer la surjectivit\'e des morphismes $j_K^*:\HH_{N+3}(Y)\to\HH_{N+1}(X_K)$
et $\bord:\HH_{N+1}(X_K)\to\HH_{N+1}(X_K,X_{A,K})$.
La premi\`ere r\'esulte de ce que $\HH_{N+1}(X_{K})^{\ev}=0$ et la
seconde de ce que $\HH_{N}(X_{A,K})_{\ev}=0$.
\cucu

\begin{lemme}\label{sd}
Le morphisme $\HH_{N-1}(X_A\setminus(X_{A,K}\setminus X_{A,K,Q}))
\to\HH_{N-1}(X_{AK+Q})$ est injectif.
\end{lemme}

\preuve
L'image par $(p_Y)_*$ du carr\'e du second diagramme du
lemme~\ref{semistable} donne le  diagramme commutatif
$$\xymatrix{
  \HH_{N-1}(X_A\setminus(X_{A,K}\setminus X_{A,K,Q}))
  \ar[r]^{\qquad F_2}\ar[d]^{v}
  &\HH_{N-1}(X_{AK+Q})\ar[d]^{(j_{AK+Q})_*}\\
  \HH_{N-1}(X_A)\ar[r]^{(j_{A})_*}&\HH_{N-1}(Y)\\},$$
donc il suffit de montrer l'injectivit\'e des morphismes  
$v:\HH_{N-1}(X_A\setminus(X_{A,K}\setminus X_{A,K,Q}))\to
\HH_{N-1}(X_A)$ et $(j_{A})_*:\HH_{N-1}(X_A)\to\HH_{N-1}(Y)$. 
La premi\`ere r\'esulte du lemme~\ref{thom}, puisque 
$\HH_{N-2}(X_{A,K},X_{A,K,Q})=0$, et la seconde de ce que
$\HH_{N-1}(X_A)_{\ev}=0$.
\cucu 

\medskip

\noindent
{\it Fin de la preuve de la
proposition~\ref{decompositionpremiere}.\ --- }  
La suite exacte qui constitue la ligne du diagramme de
la proposition~\ref{decompositionpremiere} se d\'eduit de la suite
exacte longue du lemme~\ref{semistable} et de l'isomorphisme du
lemme~\ref{contraction}. 
L'exactitude \`a gauche et \`a droite de la ligne du digramme r\'esulte des
lemmes~\ref{sg} et~\ref{sd}.
La suite exacte qui constitue la colonne
de droite est la suite exacte~(\ref{rela}).
Si $N=1$, nous avons $\HH_{N}(X_A\setminus(X_{A,K}\setminus
X_{A,K,Q})\simeq\HH_{N}(X_A\setminus(X_{A,K})$ et si $N\geq 2$
La suite exacte qui constitue la colonne de gauche est la
premi\`ere ligne du diagramme du lemme~\ref{thom}.
Le triangle de gauche du diagramme~(\ref{d1}) est le triangle du
diagramme du lemme~\ref{thom}. La commutativit\'e du triangle de
droite du diagramme~(\ref{d1}) et du 
diagramme~(\ref{d2}) s'obtiennent en prenant l'image par $(p_Y)_*$ 
respectivement du triangle du second diagramme et du premier diagramme
du lemme~\ref{semistable}.  
\cucu

\subsubsection{Lien avec les travaux de Clemens, Schmid et Steenbrink sur les 
structures de Hodge limites} \label{schmid}

La preuve de la proposition~\ref{decompositionpremiere} repose sur
l'\'etude de la famille semi-stable $\widetilde{\cal X}\to\Delta_{Q,K}$. 
Or l'\'etude de familles semi-stables du point de vue de la monodromie
autour de $\widetilde{\cal X}_0$ est l'objet d'une abondante
litt\'erature ({\it cf.}~\cite{clem} pour un langage et une situation
similaires aux notres et~\cite{sm} ou~\cite{ste} pour le cas g\'en\'eral). Notre point
de vue diff\`ere du point de vue classique en ce que nous introduisons 
au lemme~\ref{semistable} une dissym\'etrie entre les composantes de
$\widetilde{\cal X}_0$. Cette section est une digression qui compare  
(sans donner de d\'emonstration) nos r\'esultats avec les constructions
(essentiellement isomorphes) de~\cite{clem},~\cite{sm} et~\cite{ste}.

Schmid et Steenbrink d\'efinissent une structure de Hodge mixte sur la 
(co)homologie d'une fibre lisse d'une d\'eg\'enerescence semi-stable, 
qui d\'ecrit le comportement asymptotique des structures de Hodge 
lorsque la fibre tend radialement vers la fibre singuli\`ere, les 
espaces de (co)homologie des fibres \'etant identifi\'ees par 
transport plat. Dans le cas particulier de la famille $\widetilde{\cal
X}\to\Delta_{Q,K}$, cette construction donne une filtration
d\'ecroissante par le poids $W^{\bullet}$ sur $\HH_N(X_{AK+Q})$, qui
v\'erifie 
\begin{eqnarray*}
W^{-N-1}\HH_N(X_{AK+Q})&\simeq&\HH_N(X_{AK+Q});\\
W^{-N-1}/W^{-N}\HH_N(X_{AK+Q})&\simeq&\HH_{N-1}(X_{A,K})_{\ev}(1);\\
W^{-N}/W^{-N+1}\HH_N(X_{AK+Q})&\simeq&\HH_N(X_{K})^{\prim}
\oplus\HH_N(\widetilde{X_A})^{\prim};\\
W^{-N+1}\HH_N(X_{AK+Q})&\simeq&\HH_{N-1}(X_{A,K})^{\ev}.
\end{eqnarray*}
L'op\'erateur $T:\HH_N(X_{AK+Q})\to\HH_N(X_{AK+Q})$
de monodromie autour de $\widetilde{\cal X}_0$ est unipotent et
l'op\'erateur $N={\rm log}\,T=T-{\rm Id}$ est la compos\'ee de la
projection sur 
$W^{-N-1}/W^{-N}\HH_N(X_{AK+Q})$ avec un isomorphisme canonique
$$W^{-N-1}/W^{-N}\HH_N(X_{AK+Q})\simeq W^{-N+1}\HH_N(X_{AK+Q})(1).$$

La filtration $F_{\bullet}$ est une filtration de structures de Hodge
mixtes au sens de Schmid et Steenbrink. Nous avons en particulier 
\begin{eqnarray*}
W^{-N-1}/W^{-N}\HH_N(X_{AK+Q})&\simeq&F_4/F_3\HH_N(X_{AK+Q})\quad{\rm et}\\
W^{-N+1}\HH_N(X_{AK+Q})&\simeq&\HH_N(X_{A,K})^{\ev}
\stackrel{\tube}{\hookrightarrow}F_1\HH_N(X_{AK+Q}).
\end{eqnarray*}
Le gradu\'e de la filtration par le poids est donc essentiellement 
\'equivalent au gadu\'e associ\'e \`a la filtration \`a $5$ termes 
d\'eduite de la filtration $F_{\bullet}$ en coupant $X_A\setminus
X_{A,K}$ en deux morceaux. L'arbitre nous fait remarquer qu'il est
possible de d\'eduite la fitration $F_{\bullet}$ de la filtration par
le poids ci-dessus. Malheureusement, nous ne savons pas en d\'eduire
la suite exacte courte qui constitue la ligne du milieu de la
proposition~\ref{decompositionpremiere} et qui est n\'ecessaire pour
initialiser la r\'ecurrence \`a la section~\ref{preuveaux}. 

\subsubsection{La filtration quotient}\label{quot} 

Dans cette section nous supposons $A\in I_W^e$ et $Q\in I_W^d$ en plus
des hypoth\`eses de la section pr\'ec\'edente. 

Notons 
$H(W)_{A,K}$ le sous-espace de $\HH_{N-2}(X_{A,K,Q})$ engendr\'e par les
classes des composantes irr\'eductibles de $W\cap X_K$ et par 
$(j_{A,K,Q}^{A,K})^*\HH_N(X_{A,K})$. Les espaces 
$H(W)\subset\HH_N(X_{AK+Q})$ et $H(W)_{A,K}\subset\HH_{N-2}(X_{A,K,Q})$
sont $\pi_1({\cal U}_{A,Q},{\bf K})$-invariants. 

\begin{lemme}\label{qhw}
Si $\hh_{N-1}(X_{A,K})_{\ev}\neq 0$ alors
la filtration $F_{\bullet}\HH_N(X_{AK+Q})$ de 
$\pi_1({\cal U}_{A,Q},K)$-modules donne par restriction \`a
$H(W)^{\perp}$ la filtration suivante.
\begin{enumerate}
\item $F_1H(W)^{\perp}\simeq\HH_{N+1}(X_{A}\setminus X_{A,K})$,
\item $F_2/F_1H(W)^{\perp}\simeq H(W)_{A,K}^{\perp}$,
\item $F_3/F_2H(W)^{\perp}\simeq\HH_{N+2}(X_{K})^{\ev}$ ,
\item $F_4/F_3H(W)^{\perp}\simeq\HH_{N+1}(X_{A,K})_{\ev}$, 
\end{enumerate}
et $F_4H(W)^{\perp}=H(W)^{\perp}$. 
\end{lemme}

Gr\^ace \`a la suite exacte~(\ref{ouvert}) (pour l'assertion~1),
nous  obtenons imm\'ediatement~:

\begin{cor}\label{qhw1} Nous avons
\begin{enumerate}
\item $\dim F_1 H(W)_j^{\perp}=
      \hh_{N-2j+2}(X_{A,j})^{\prim}+\hh_{N-2j+1}(X_{A,K,j})_{\ev}$~;
\item $\dim F_2/F_1H(W)_j^{\perp}\leq\hh_{N-2j}(X_{A,K,Q,j})_{\ev}$~;
\item $\dim F_3/F_2H(W)_j^{\perp}=\hh_{N-2j+2}(X_{K,j})_{\ev}$~;
\item $\dim F_4/F_3H(W)_j^{\perp}=\hh_{N-2j+1}(X_{A,K,j})_{\ev}$.
\end{enumerate}
\end{cor}

\preuve 
Nous \'etudions l'image de $H(W)$ par la filtration $F_{\bullet}$.
Notons $[W]\subset\HH_N(X_{A,Q})$ l'espace vectoriel engendr\'e par les
classes des composantes irr\'eductibles de $W$ et posons
$$H(W)'=(j_{A,Q}^{AK+Q})_*[W]+
 j_{AK+Q}^*\circ c_1(\anneau(1))\frown\HH_{N+4}(Y).$$
Comme le morphisme 
$j_{AK+Q}^*\circ c_1(\anneau(1))\frown:\HH_{N+4}(Y)\to\HH_N(X_{AK+Q})$ 
s'identifie \`a la multiplication par un scalaire pr\`es au morphisme
$(j_{A,Q}^{AK+Q})_*\circ (j_{A,Q})^*$, 
nous avons  $H(W)'\subset(j_{A,Q}^{AK+Q})_*\HH_N(X_{A,Q})$.
Comme le triangle de droite du diagramme~(\ref{d1}) est commutatif, nous
avons $(j_{A,Q}^{AK+Q})_*=F_2\circ\phi$, donc
$H(W)'\subset F_2\HH_N(X_{AK+Q})$, et, comme
$\hh_{N-1}(X_{A,K})_{\ev}\neq 0$, d'apr\`es le lemme~\ref{thom} nous
avons $\im(\alpha)\cap\im(\phi)=0$, donc
$H(W)'\cap F_1\HH^N(X_{AK+Q})=0$.
La projection de $H(W)'$ sur $\HH_{N-2}(X_{A,K,Q})$
est donc injective. Comme le triangle de gauche du
diagramme~(\ref{d1}) est commutatif et que  
$j_{A,K,Q}^*\HH_{N+4}(Y)=(j_{A,K,Q}^{A,K})^*\HH_N(X_{A,K})$,
nous avons
\begin{eqnarray}\label{bout1}
H(W)'\simeq H(W)_{A,K}\subset F_2/F_1\HH_N(X_{AK+Q}).
\end{eqnarray}

D'autre part, nous consid\'erons l'espace
$$j_{AK+Q}^*\HH_{N+2}(Y)^{\prim}=
  \coker(j_{AK+Q}^*\circ c_1(\anneau(1))\frown\HH_{N+4}(Y)\to 
  j_{AK+Q}^*\HH_{N+2}(Y)),$$ 
o\`u la fl\`eche est l'inclusion. Le morphisme $j_K^*$ d\'efinit 
par passage au quotient un isomorphisme  
$\HH_{N+2}(Y)^{\prim}\simeq
  \kker\left(\HH_{N}(X_K)^{\prim}\to\HH_{N}(X_K)^{\ev}\right)$.
Les morphismes $j_{AK+Q}^*:\HH_{N+2}(Y)\to\HH_{N}(X_{AK+Q})$ et 
$\rel:\HH_{N}(X_K)^{\prim}\to\HH_{N}(X_K,X_{AK})$ sont injectifs.
Comme le diagramme~(\ref{d2}) est commutatif, nous avons donc
\begin{eqnarray}\label{bout2}
j_{AK+Q}^*\HH_{N+2}(Y)^{\prim}\simeq
\kker\left(\HH_{N}(X_K)^{\prim}\to\HH_{N}(X_K)^{\ev}\right)\subset
F_3/F_2\HH_N(X_{AK+Q}).
\end{eqnarray}

Enfin, comme $(j_{A,Q}^{AK+Q})_*[W]$ est inclus dans
$F_2\HH^N(X_{AK+Q})$ 
et que $j_{AK+Q}^*\HH_{N+2}(Y)^{\prim}$ s'injecte dans 
$F_3/F_2\HH_N(X_{AK+Q})$, de la suite exacte
$$0\to j_{AK+Q}^*\circ c_1(\anneau(1))\frown\HH_{N+4}(Y)\to 
  j_{AK+Q}^*\HH_{N+2}(Y)\to j_{AK+Q}^*\HH_{N+2}(Y)^{\prim}\to 0$$
nous d\'eduisons donc la suite exacte
\begin{eqnarray*}\label{lesdeuxbouts}
0\to H(W)'\to H(W)\to j_{AK+Q}^*\HH_{N+2}(Y)^{\prim}\to 0,
\end{eqnarray*}
et le lemme~\ref{qhw} r\'esulte des assertions~(\ref{bout1})
et~(\ref{bout2}).   
\cucu 

\section{Preuve du th\'eor\`eme~\ref{mono}}

\subsection{Preuve du th\'eor\`eme~\ref{mono}}

Nous adoptons les notations de la section~\ref{intro} et faisons les
hypoth\`eses du th\'eor\`eme~\ref{mono}~; nous supposons de plus 
${\cal V}^{d}(W)$ non vide. 

Pour tout $F\in{\cal V}^{d}(W)$ nous identifions $\HH_{N}(X_F)$ \`a
$\HH^{N}(X_F)$ par la dualit\'e de Poincar\'e. 
Pour tout $\lambda\in H(W)^{\perp}$, nous notons $H(\lambda)$ la
$\pi_1({\cal V}^{d}(W),F)$-repr\'esentation de monodromie engendr\'ee par $\lambda$,
$\hh(\lambda)$ sa dimension et $\hh(\lambda)^{\perp}$ sa codimension 
dans $H(W)^{\perp}$.  

\medskip

L'argument central de la preuve consiste \`a minorer $\hh(\lambda)$ pour
tout $\lambda\in H(W)^{\perp}\setminus\{0\}$. C'est l'objet de la
proposition suivante qui synth\'etise les r\'esultats des deux premi\`eres
parties et dont la preuve occupe la section~\ref{preuveaux}.   

\begin{prop} \label{aux}
Pour tout $F\in{\cal V}^{d}(W)$ et pour tout
$\lambda\in H(W)^{\perp}\setminus\{0\}$ nous avons 
\begin{enumerate}
\item 
$$\hh(\lambda)\geq {\rm min}\left(\left.\begin{array}{c}
\hh_{N-2i+2}(Y(\underbrace{e,\dots,e}_{i-1\ {\rm fois}},
               \underbrace{d-e,\dots,d-e}_{i\ {\rm fois}})_{\ev}\\
\hh_{N-2i+1}(Y(\underbrace{e,\dots,e}_{i\ {\rm fois}},
               \underbrace{d-e,\dots,d-e}_{i\ {\rm fois}})_{\ev}
		 \end{array}\right|\ 
i\in\left\{1,\dots, \left[\frac{N+1}{2}\right]\right\}\right),$$
o\`u $Y(d_1,\dots d_r)$ d\'esigne une intersection compl\`ete lisse de
multidegr\'e $(d_1,\dots d_r)$~;
\item si $\hh(\lambda)>A_{Y,e}$ alors $\hh(\lambda)^{\perp}\leq A_{Y,e}$, 
o\`u nous avons pos\'e
$$A_{Y,e}=\hh_N(Y(e))^{\prim}+2\hh_{N-1}(Y(e,d-e))_{\ev}+
          \hh_{N-2}(Y(d,e,d-e))_{\ev}.$$
\end{enumerate}
\end{prop}

La proposition suivante donne une estimation de la dimension des
espaces d'homologie intervenant dans la proposition~\ref{aux}. Sa 
preuve occupe la section~\ref{preuveestimation}. 

\begin{prop}\label{estimation}
Il existe une constante $C\in\RR_+^*$ qui ne d\'epend que de $Y$ et de 
$\anneau(1)$, telle que pour tout $r\in\{1,\dots,N+1\}$, pour tout
$r$-uplet d'entiers strictement positifs et croissants
$(d_1,\dots,d_r)$ tels que $d_r\geq C$ et pour toute intersection
compl\`ete lisse $X\subset Y$ de multidegr\'e $(d_1,\dots ,d_r)$ nous avons
\begin{eqnarray*}
\frac{1}{2}d_r^k\prod_{i=1}^rd_i\leq\hh_k(X)_{\ev}
\leq\hh_k(X)\leq Cd_r^k\prod_{i=1}^rd_i.
\end{eqnarray*}
o\`u nous avons pos\'e $k=\dim X=N+1-r$.
\end{prop}

\noindent
{\it Fin de la preuve du th\'eor\`eme~\ref{mono}.\ ---\ }
Comme la $\pi_1({\cal V}^{d}(W),F)$-repr\'esentation de monodromie sur $H(W)^{\perp}$ 
pr\'eserve la forme d'intersection, elle est somme directe orthogonale
de ses sous-repr\'esentations. Elle est donc irr\'eductible si
et seulement si pour tout $\lambda$ et $\mu$ dans
$H(W)^{\perp}\setminus\{0\}$ nous avons $H(\mu)\cap H(\lambda)\neq0$.

Choisissons $\lambda$ et $\mu$ dans $H(W)^{\perp}\setminus\{0\}$
et supposons $d\geq2e$. 

D'apr\`es l'assertion 1 de la proposition~\ref{aux}
(pour $e={\left[\frac{d}{2}\right]}$), $\hh(\lambda)$ est minor\'e par
la dimension de l'homologie en dimension moiti\'e d'une intersection
compl\`ete lisse d'hypersurfaces de degr\'es
${\left[\frac{d}{2}\right]}$ et $d-{\left[\frac{d}{2}\right]}$~; comme
$d-{\left[\frac{d}{2}\right]}\geq{\left[\frac{d}{2}\right]}$,
d'apr\`es la proposition~\ref{estimation}, il
existe une constante $C_1$ qui ne d\'epend que de $Y$ et de
$\anneau(1)$ telle que, si $d\geq C_1$, nous avons  
$\hh(\lambda)\geq\frac{1}{2}{\left[\frac{d}{2}\right]}^{N+1}$ et
$\hh(\mu)\geq\frac{1}{2}{\left[\frac{d}{2}\right]}^{N+1}$. 

D'autre part, d'apr\`es la proposition~\ref{estimation} pour les
intersections compl\`etes de multidegr\'es $(e)$, $(e,d-e)$ et
$(d,e,d-e)$, il existe une constante $C_2$ qui ne d\'epend que de $Y$ et
de $\anneau(1)$ telle que $A_{Y,e}\leq C_2ed^{N}$. Nous avons donc
$\frac{1}{2}{\left[\frac{d}{2}\right]}^{N+1}>C_2ed^{N}\geq A_{Y,e}$
pour $d\geq 3^{N+2}C_2e$.

Supposons $d\geq {\rm max}\,(C_1,3^{N+2}C_2)e$~: nous avons
$\hh(\lambda)> A_{Y,e}$, $\hh(\mu)> A_{Y,e}$ et, d'apr\`es
l'assertion~2 de la proposition~\ref{aux},
$\hh(\mu)^{\perp}\leq A_{Y,e}$, donc $\hh(\lambda)>\hh(\mu)^{\perp}$
et $H(\mu)\cap H(\lambda)\neq0$.
\cucu

\subsection{Preuve de la proposition~\ref{aux}}
\label{preuveaux}

Comme ${\cal V}^{d}(W)$ est connexe, il est suffisant de montrer la
proposition~\ref{aux} pour un seul $F\in{\cal V}^{d}(W)$ et, au lieu de
${\cal V}^{d}(W)$, pour un sous-espace de ${\cal V}^{d}(W)$ contenant $F$.
L'assertion~1 de la proposition~\ref{aux} est donc impliqu\'ee par les 
lemmes~\ref{lissite} et~\ref{recu} \'enonc\'es \`a la
section~\ref{secrecu}. La preuve par r\'ecurrence du lemme~\ref{recu}
occupe les sections~\ref{secrecu1} et~\ref{secrecu2}. 

La preuve de la seconde assertion ne n\'ecessite pas de r\'ecurrence~;
elle occupe la section~\ref{secrecu3}.

\subsubsection{Preuve de l'assertion~1~: l'\'enonc\'e de la r\'ecurrence}
\label{secrecu}

Posons $r=\left[\frac{N+1}{2}\right]$.

Comme le sch\'ema $W\subset Y$ v\'erifie l'assertion $(a)$ de la
proposition~\ref{bertini}, les assertions~$(c)$ et~$(d)$ de la
proposition~\ref{bertini} sont vraies~: pour
$(A_1,\dots,A_r,R)\in(I_W^e)^r\times I_W^d$ g\'en\'erique 
\begin{enumerate}
\item $A_1,\dots,A_r$ et $R$ sont en intersection compl\`ete~;
\item pour tout $j\in\{1,\dots,r\}$ le lieu singulier des sch\'emas
$$X_{A_1}\cap\cdots\cap X_{A_j}\quad{\rm et}\quad 
  X_{A_1}\cap\cdots\cap X_{A_j}\cap X_R$$ 
  est support\'e par $W$ et est de
  dimension inf\'erieure ou \'egale \`a $j-2$ et $j-1$ respectivement
  (par convention, un lieu de dimension $-1$ est vide)~;
\item le sch\'ema $(X_{A_1}\cap\cdots\cap X_{A_r}\cap X_R)\setminus W$ 
  est lisse et connexe, sauf si $N=1$.
\end{enumerate}
Nous fixons $(A_1,\dots,A_r,R)\in(I_W^e)^r\times I_W^d$ comme
ci-dessus.

\medskip

Soit $j\in\{1,\dots,r\}$.
Pour tout $(K_1,\dots,K_{j-1})\in\HH^0(Y,(\anneau(d-e))^{j-1}$
posons 
$$Y_j=\bigcap_{i=1}^{j-1} X_{A_i,K_i},\quad{\rm et}\quad
W_j=W\cap\bigcap_{i=1}^{j-1} X_{K_i}.$$
Pour tout $i\in\NN$ notons 
\begin{eqnarray*}
\psi^{i}_{j}:\HH^0(Y,(\anneau(i))&\to&\HH^0(Y_j,\anneau(i))
\quad {\rm et}\quad\\
{}^A\psi^{i}_{j}:\HH^0(Y,(\anneau(i))&\to&
                 \HH^0(Y_j\cap X_{A_{j}},\anneau(i))
\end{eqnarray*}
les applications de restriction~; notons
${\cal V}_j^i\subset\HH^0(Y_j,\anneau(i))$ et
${}^A{\cal V}_j^i\subset\HH^0(Y_j\cap X_{A_{j}},\anneau(i))$
les ouverts param\'etrant (\`a la multiplication par un scalaire pr\`es)
les hypersurfaces lisses de $Y_j$ et de $Y_j\cap X_{A_{j}}$
de classe $c_1(\anneau (i))$. 
Pour tout $P\in\HH^0(Y_j,\anneau(i))$ notons $X_{P,j}\subset Y_j$
l'hypersurface associ\'ee.

Pour tout $(K_j,\dots,K_r)\in\HH^0(Y,(\anneau(d-e))^{r-j+1}$ posons
$${\cal W}_j^d=\left\{\left. F\in{\cal V}_j^d\ \right|\
F=\psi^{d}_{j}\left(\sum_{i=j}^rA_iK_i+R\right),\ 
(K_j,\dots,K_r)\in\HH^0(Y,(\anneau(d-e))^{r-j+1}\right\}.$$
Remarquons que ${\cal W}_1^d\subset{\cal V}^d(W)$.

Nous dirons que $(K_1,\dots,K_{j-1})$ est {\em rouge}
si les vari\'et\'es $Y_j$ et $Y_j\cap X_{A_{j}}$ sont lisses, 
de dimensions respectivement $N-2j+3$ et $N-2j+2$,
et si l'espace ${\cal W}_j^d$ est non vide. 

\begin{lemme}\label{lissite}
Il existe un $(K_1,\dots,K_{j-1})\in\HH^0(Y,(\anneau(d-e))^{j-1}$ rouge.
\end{lemme}

\preuve 
D'apr\`es le th\'eor\`eme de Bertini usuel et les
conditions 1 et 2 sur les $A_i$ et $R$, pour $(K_1,\dots,K_{j-1})$
g\'en\'erique les vari\'et\'es $Y_j$, $Y_j\cap X_{A_j}$ et 
$Y_j\cap X_R$ sont lisses, de dimensions respectives $N-2j+3$,
$N-2j+2$ et $N-2j+2$. En particulier, nous avons 
$\psi^{d}_{j}(R)\in{\cal W}_j^d$. Donc $(K_1,\dots,K_{j-1})$ g\'en\'erique
est rouge. 
\cucu

\medskip

Pour tout $(K_1,\dots,K_{j-1})$ rouge et
$F\in{\cal W}_j^d$ notons $H(W)_j\subset\HH_{N-2j+2}(X_{F,j})$ le 
sous-espace engendr\'e par les composantes irr\'eductibles de $W_j$ et
par l'image de $\HH_{N-2j+4}(Y_{j})$. Le groupe $\pi_1({\cal W}_j^d,F)$
agit alors par monodromie sur $\HH_{N-2j+2}(X_{F,j})$. 

Nous pouvons enfin donner l'\'enonc\'e de la r\'ecurrence~: 

\begin{lemme}\label{recu}
Pour tout $j\in\{1,\dots,r\}$, $(K_1,\dots,K_{j-1})$ rouge et 
$F\in{\cal W}_j^d$, les dimensions des 
$\pi_1({\cal W}_j^d,F)$-repr\'esentations
de monodromie engendr\'ees par les \'el\'ements de 
$H(W)_j^{\perp}\setminus\{0\}$ v\'erifient la minoration de
l'assertion~1.
\end{lemme}

\subsubsection{Le lemme~\ref{recu} pour $j+1$ implique le 
lemme~\ref{recu} pour $j$}
\label{secrecu1}

Nous fixons $j\in\{1,\dots,r\}$ et $(K_1,\dots,K_{j-1})$ rouge.

Pour tout $(K_j,\dots,K_r)\in\HH^0(Y,(\anneau(d-e))^{r-j+1}$ posons
$$A=\psi^{-e}_{j}(A_j),\quad K=\psi^{d-e}_{j}(K_j),\quad 
  Q=\psi^{d}_{j}\left(\sum_{i=j+1}^rA_iK_i+R\right)\quad {\rm et}\quad 
  F=AK+Q.$$ 

Nous dirons que $(K_j,\dots,K_r)\in\HH^0(Y,(\anneau(d-e))^{r-j+1}$ est 
{\em bleu} si les vari\'et\'es 
$$X_{A,j},\  X_{K,j},\  X_{F,j},\  X_{A,K,j}=Y_{j+1},\  
X_{A,K,Q,j}=Y_{j+1}\cap X_{Q,j}\ {\rm et}\ Y_{j+1}\cap X_{A_{j+1}}$$ 
sont lisses, de codimensions attendues.
Remarquons que les trois derni\`eres conditions impliquent que 
$(K_1,\dots,K_{j-1},K_j)$ est rouge.

Pour tout $(K_{j+1},\dots,K_r)\in\HH^0(Y,(\anneau(d-e))^{r-j}$
notons ${\cal V}_{\rm bleu}^{d-e}\subset\HH^0(Y,(\anneau(d-e))$
l'espace des $K_j$ tels que $(K_j,K_{j+1},\dots,K_r)$ est bleu.
D'apr\`es le th\'eor\`eme de Bertini usuel, 
$(K_j.K_{j+1},\dots,K_r)$ g\'en\'erique
est bleu~; fixons
$(K_{j+1},\dots,K_r)\in\HH^0(Y,(\anneau(d-e))^{r-j}$
tel que lequel l'espace ${\cal V}_{\rm bleu}^{d-e}$ est non vide.

Pour tout $K_j\in{\cal V}_{\rm bleu}^{d-e}$
les vari\'et\'es indic\'ees par $j$ v\'erifient alors les conditions
de la situation g\'en\'erale de la section~\ref{referee}. 
Nous pouvons donc raisonner comme dans la section~\ref{rec}~:
la vari\'et\'e $X_{F,j}$ est lisse, donc pour 
$t\in\PP_{\C}^1$ g\'en\'erique la vari\'et\'e $X_{AK+tQ,j}$ est
lisse. Soit $\rho_{K}$ le plus petit r\'eel 
(\'eventuellement infini) tel que
si $\Delta\subset\C$ est un disque ouvert de rayon $\rho_{K}$ et de
centre~$0$, alors pour tout $t\in\Delta\setminus\{0\}$ et
$j\in\{1,\dots,r\}$ la vari\'et\'e $X_{AK+tQ,j}$ est lisse.
Posons
$${\cal U}=\left\{\left.K_j\in{\cal V}_{\rm bleu}^{d-e}\
  \right|\ \rho_{K}>1\right\}.$$

Le groupe $\pi_1({\cal U},K_j)$ agit par monodromie sur 
$\HH_{N-2j+2}(X_{F,j})$
via l'application ${\cal U}\to{\cal W}_j$, $K_j\mapsto F_j$. 
Notons ${\cal V}_{A,Q,j}^{d-e}$, ${}^A{\cal V}_{j}^{d-e}$,
${}^A{\cal V}_{Q,j}^{d-e}$ et ${\cal U}_{A,Q,j}$
les analogues
des espaces ${\cal V}_{A,Q}^{d-e}$, ${}^A{\cal V}^{d-e}$,
${}^A{\cal V}_{Q}^{d-e}$ et ${\cal U}_{A,Q}$
des sections~\ref{referee} et~\ref{rec} pour les
vari\'et\'es indic\'ees par $j$. 
Nous avons $\psi^{e}_{j}({\cal U})\subset{\cal U}_{A,Q,j}$, 
donc le groupe $\pi_1({\cal U},K_j)$ agit sur les espaces 
\begin{itemize}
\item $\HH_{N-2j+2}(X_{A,j}\setminus (X_{A,K,j}\setminus X_{A,K,Q,j}))$
      via $\pi_1({}^A{\cal V}_{Q,j}^{d-e},K)$,
\item $\HH_{N-2j+2}(X_{A,j}\setminus X_{A,K,j})$
      via $\pi_1({}^A{\cal V}_{j}^{d-e},K)$,
\item $\HH_{N-2j+2}(X_{K,j})^{\ev}$ 
      via $\pi_1({\cal V}_j^{d-e},K)$ et
\item $\HH_{N-2j+1}(X_{A,K,j})_{\ev}$ 
      via $\pi_1({}^A{\cal V}_{j}^{d-e},K)$.
\end{itemize}

\begin{lemme} \label{existance}
Les applications $\psi^{e}_{j}:{\cal U}\to{\cal V}_j^{d-e}$,
${}^A\psi^{e}_{j}:{\cal U}\to{}^A{\cal V}_{Q,j}^{d-e}$ et 
${}^A\psi^{e}_{j}:{\cal U}\to{}^A{\cal V}_j^{d-e}$
induisent des morphismes de groupes fondamentaux surjectifs.
\end{lemme}

\preuve Le lemme r\'esulte des faits suivants~:
\begin{itemize}
\item  
$\C^*{\cal U}={\cal V}_{\rm bleu}^{d-e}$, donc ${\cal U}$ et
${\cal V}_{\rm bleu}^{d-e}$ sont homotopes~; 
\item 
les applications $\psi^{e}_{j}:{\cal V}_{\rm bleu}^{d-e}\to
{\cal V}_{j}^{d-e}\cap\im\psi^{e}_{j}$ et
${}^A\psi^{e}_{j}:{\cal V}_{\rm bleu}^{d-e}\to
{}^A{\cal V}_{Q,j}^{d-e}\cap\im{}^A\psi^{e}_{j}$ sont dominantes,
\`a fibres connexes~;
\item 
les inclusions
${\cal V}_{j}^{d-e}\cap\im\psi^{e}_{j}\subset{\cal V}_{j}^{d-e}$ et
${}^A{\cal V}_{Q,j}^{d-e}\cap\im\psi^{e}_{j}\subset{}^A{\cal V}_{Q,j}^{d-e}$ 
induisent des morphismes de groupes fondamentaux surjectifs par
le th\'eor\`eme de Zariski~: en effet, les applications 
$\psi^{e}_{j}$ et ${}^A\psi^{e}_{j}$ ne sont en g\'en\'eral pas 
surjectives, mais comme $\anneau(d-e)$ est tr\`es ample, les images 
de $\psi^{e}_{j}$ et ${}^A\psi^{e}_{j}$ s'identifient \`a l'espace des
sections hyperplanes des vari\'et\'es $Y_j$ et $X_{A,j}$ pour 
le plongement projectif de $Y$ d\'efini par $\anneau(d-e)$~; 
en particulier, les images de $\psi^{e}_{j}$ et ${}^A\psi^{e}_{j}$
rencontrent transversalement les hypersurfaces
discriminantes des espaces $\HH^0(Y_j,\anneau(d-e))$ et
$\HH^0(X_{A,j},\anneau(d-e))$, o\`u les vari\'et\'es
$X_{K,j}$, $X_{A,K,j}$ et $X_{A,K,Q,j}$ sont singuli\`eres~;
\item 
${}^A{\cal V}_{Q,j}^{d-e}\subset{}^A{\cal V}_j^{d-e}$ est
un ouvert de Zariski.\cucu
\end{itemize}

\medskip

\noindent
Pour toute $\pi_1({\cal W}_j^d,F)$- et 
$\pi_1({\cal U},K_j)$-repr\'esentation $E$ notons ${\rm m}^{\cal W}_jE$ 
et ${\rm m}_jE$ respectivement
la dimension minimale d'une sous-repr\'e\-sen\-ta\-tion 
non nulle (si $E=0$ nous convenons que ${\rm m}_jE=\infty$). Nous avons

\begin{lemme}\label{qhw2} Si $\hh_{N-2j+1}(X_{A,K,j})_{\ev}\neq 0$ et 
$\hh_{N-2j+2}(X_{K,j})_{\ev}\neq 0$ alors
\begin{enumerate}
\item ${\rm m}_j\,F_1H(W)_j^{\perp}\geq\hh_{N-2j+1}(X_{A,K,j})_{\ev}$~;
\item ${\rm m}_j\,F_2/F_1H(W)_j^{\perp}=
       {\rm m}^{\cal W}_{j+1}\,H(W)_{j+1}^{\perp}$~;
\item ${\rm m}_j\,F_3/F_2H(W)_j^{\perp}=\hh_{N-2j+2}(X_{K,j})_{\ev}$~;
\item ${\rm m}_j\,F_4/F_3H(W)_j^{\perp}=\hh_{N-2j+1}(X_{A,K,j})_{\ev}$.
\end{enumerate}
\end{lemme}

\preuve D'apr\`es les lemmes~\ref{qhw} et~\ref{existance},
\begin{enumerate} 
\item r\'esulte de la proposition~\ref{monodromie-sur-complementaire}~; 
\item r\'esulte de ce que, comme $(K_j,\dots,K_r)$ est bleu, 
      $(K_1,\dots,K_{j-1},K_j)$ est rouge et 
      $X_{A,K,Q,j}\subset Y_{j+1}$ est l'hypersurface
      d\'efinie par 
      $\psi^{d}_{j+1}\left(\sum_{i=j+1}^rA_iK_i+R\right)\in{\cal W}_{j+1}^d$~;
\item r\'esulte de l'irr\'eductibilit\'e de l'action de
      $\pi_1({\cal V}_j^{d-e},K)$ sur $\HH_{N-2j+2}(X_{K,j})^{\ev}$~; 
\item r\'esulte de l'irr\'eductibilit\'e de l'action de
      $\pi_1({\cal V}_j^{d-e},K)$ sur $\HH_{N-2j+1}(X_{A,K,j})^{\ev}$.\cucu
\end{enumerate}

\medskip

\noindent
Nous devons minorer ${\rm m}^{\cal W}_jH(W)_j^{\perp}$. Comme
${\rm m}^{\cal W}_jH(W)_j^{\perp}\geq {\rm m}_jH(W)_j^{\perp}$,
le lemme~\ref{existance} implique l'\'enonc\'e suivant, 
qui montre que le lemme~\ref{recu} pour $j+1$ implique 
le lemme~\ref{recu} pour $j$~:

\begin{cor}\label{dimdim}
Pour tout $j\in\{1,\dots,r\}$, 
$${\rm m}^{\cal W}_j\,H(W)_j^{\perp}\geq{\rm min}\left(
\hh_{N-2j+1}(X_{A,K,j})^{\ev},\ 
\hh_{N-2j+2}(X_{K,j})^{\ev},\
{\rm m}^{\cal W}_{j+1}\,H(W)^{\perp}_{j+1}\right).$$
\end{cor}

\subsubsection{Preuve du lemme~\ref{recu} pour $j=r$}
\label{secrecu2}

\noindent
{\it Cas $N=2r-1$.\ ---\ } Nous avons $H(W)^{\perp}_{r+1}=0$
et le corollaire du lemme~\ref{qhw2} implique l'assertion~1.
\cucu

\medskip

\noindent
{\it Cas $N=2r$.\ ---\ } Nous ne pouvons appliquer le
corollaire du lemme~\ref{qhw2} parce que 
${\rm m}^{\cal W}_{j+1}\,H(W)^{\perp}_{r+1}=
{\rm m}_j\,F^1/F^2H(W)^{\perp}_{r}$ est
trop petit. Le corollaire du lemme~\ref{qhw2} implique seulement qu'il 
suffit de minorer ${\rm m}_j\, F^2H(W)^{\perp}_{r}$.

\begin{lemme}\label{N=2H(W)}
L'espace $H(W)_r\subset\HH_2(X_{F,r})$ contient les classes de 
toutes les composantes irr\'eductibles de la courbe $X_{A,Q,r}$. 
\end{lemme}

\preuve
Le sch\'ema $W_r$ est un sous-sch\'ema de la courbe
$X_{A,Q,r}$. Notons $W_r'$ l'adh\'erence de $X_{A,Q,r}\setminus W_r$
dans $X_{A,Q,r}$~: nous avons  $X_{A,Q,r}=W_r\cup W_r'$. Par
la condition $3$ sur les $A_i$, le sch\'ema $W_r'\setminus(W_r'\cap
W_r)$ est connexe et lisse, donc le
sch\'ema $W_r'$ est une courbe irr\'eductible. D'une part les classes
des composantes irr\'eductibles de $W_r$ appartiennent par d\'efinition
\`a $H(W)_r$, d'autre part, comme $X_{A,Q,r}=X_{A,r}\cap X_{F,r}$, la
classe de $X_{A,Q,r}$ est la  restriction \`a $H_2(X_{F,r})$ de la
classe de $X_{A,r}$ dans $H_4(Y_r)$, donc appartient aussi \`a $H(W)_r$. 
Ceci implique que la classe de $W_r'$ appartient \`a $H(W)_r$. 
\cucu

\medskip

L'espace $F_2\HH_2(X_{Q,r})$ s'identifie \`a 
$\HH_2(X_{A,r}\setminus(X_{A,K,r}\setminus X_{A,K,Q,r}))$
par la proposition~\ref{decompositionpremiere} pour les vari\'et\'es
indic\'ees par $r$. Comme le triangle de droite du
diagramme~(\ref{d1}) est commutatif, 
le lemme~\ref{N=2H(W)} implique
$\im\phi\subset F_2H(W)_r$, donc soit $F_2H(W)_r^{\perp}=0$ soit
$F_2H(W)_r^{\perp}$ rencontre 
$\HH_N(X_{A,r}\setminus(X_{A,K,r}\setminus X_{A,K,Q,r}))\setminus\im\phi$. 
La proposition~\ref{surfaces} et le lemme~\ref{existance} impliquent alors 
${\rm m}_j\,F^2H(W)_r^{\perp}\geq\hh_{1}(X_{A,K,r})_{\ev}$
et le lemme~\ref{recu} est vrai.   
\cucu

\subsubsection{Preuve de l'assertion~2}
\label{secrecu3}

D'apr\`es le corollaire du lemme~\ref{qhw}, nous avons 
\begin{eqnarray*} 
A_{Y,e}&\geq&\dim F_1 H(W)^{\perp}+\dim F_2/F_1H(W)^{\perp}
                                    +\dim F_4/F_3H(W)^{\perp}\\
&\geq&\dim H(W)^{\perp}-\dim F_3/F_2H(W)^{\perp}.
\end{eqnarray*}
La condition $\hh(\lambda)>A_{Y,e}$ implique donc 
$F_3/F_2H(\lambda)\neq 0$. Or d'apr\`es le lemme~\ref{qhw2}, 
${\rm m}_1F_3/F_2H(W)^{\perp}=\dim F_3/F_2H(W)^{\perp}$, donc
$F_3/F_2H(\lambda)=F_3/F_2H(W)^{\perp}$, donc
$\hh(\lambda)\geq \dim F_3/F_2H(W)^{\perp}$ et
$\hh(\lambda)^{\perp}\leq \dim H(W)^{\perp}-\dim F_3/F_2H(W)^{\perp}
                 \leq A_{Y,e}.$
\cucu

\subsection{Preuve de la proposition~\ref{estimation}~: les estimations
asymptotiques} \label{preuveestimation} 

Soit $\chi(X)=\sum_{i=0}^{2k}(-1)^i\hh_i(X)$ la
caract\'eristique d'Euler de $X$, soient $T_X$ et $T_Y$ les fibr\'es
tangents respectifs de $X$ et de $Y$, et soient 
$c\,(T_X)$ et $c\,(T_Y)$ leurs classes de Chern totales. 
Nous posons $\eta=c_1(\anneau(1))$ et 
$c\,(T_Y)=\sum_{i=0}^{N+1}c_i(T_Y)$, et nous notons
$Q_i(d_1,\dots,d_r)$ le coefficient de degr\'e $k-i$ en $t$ de la
s\'erie enti\`ere $\frac{1}{(1+d_1t)\dots (1+d_rt)}$. Nous avons alors
\begin{eqnarray*}
\chi(X)&=&\int_X c\,(T_X)
=\int_X\frac{c\,(T_Y)}{(1+d_1\eta)\dots (1+d_r\eta)}
=\prod_{i=1}^rd_i\int_Y  
\frac{\eta^rc\,(T_Y)}{(1+d_1\eta)\dots (1+d_r\eta)}\\
      &=&\left(\prod_{i=1}^rd_i\right)
  \sum_{i=0}^kQ_i(d_1,\dots,d_r)\int_Y\eta^{N+1-i}c_i(T_Y).
\end{eqnarray*}
D'une part il existe une constante $C_r\in \RR_+^*$ qui ne d\'epend que
de $r$ telle que pour tout $i\in\{0,\dots,k\}$ nous avons
$d_r^{k-i}\leq |Q_i(d_1,\dots,d_r)|\leq C_rd_r^{k-i}.$
D'autre part, il existe une constante
$C_T\in \RR_+^*$  qui d\'epend de $Y$ et de $\anneau_Y(1)$ telle que
pour tout $i\in\{1,\dots,r\}$ nous avons 
$\int_Y\eta^{N+1-i}c_i(T_Y)\leq C_T$.
Posons $C_{\chi}={\rm max}\,((N+1)C_TC_r\ |\ r\in\{1,\dots,N+1\}).$
Comme $\int_Y\eta^{N+1}c_0(T_Y)=\int_Y\eta^{N+1}$, nous
avons $1\leq\int_Y\eta^{N+1}c_0(T_Y)$. Nous en d\'eduisons
\begin{eqnarray}\label{minmaj} 
\left(1-\frac{C_{\chi}}{d_r}\right)d_r^k\prod_{i=1}^rd_i\leq |\chi(X)|
\leq C_{\chi}d_r^k\prod_{i=1}^rd_i.
\end{eqnarray}
 
Or la th\'eorie de Lefschetz donne
\begin{eqnarray*}
\hh_i(X)&=&\hh_i(Y) \qquad {\rm pour\ }i\in\{0,\dots,k-1\}\ ; \\
\hh_k(X)&=&\hh_k(X)_{\ev}+\hh_k(Y)\ ; \\
\hh_i(X)&=&\hh_{2k-i}(X)=\hh_{2k-i}(Y)=\hh_{2r+i}(Y) 
          \qquad {\rm pour\ }i\in\{k+1,\dots,2k\}.
\end{eqnarray*}
Posons $C_H=\sum_{i=0}^{2N+2}\hh_i(Y)$. Nous avons donc
$$(-1)^k\chi(X)-C_H\leq\hh_k(X)_{\ev}\leq\hh_k(X)\leq(-1)^k\chi(X)+C_H.$$
La proposition~\ref{estimation} r\'esulte de cette in\'egalit\'e
et de l'in\'egalit\'e~(\ref{minmaj}) pour $d_r\geq 2C_{\chi}+2C_H$. 
\cucu

\medskip

\noindent
Ania Otwinowska\\
Laboratoire de Math\'ematiques\\
Universit\'e Paris-Sud - B\^at 425\\
91405 Orsay Cedex\\
France\\
e-mail~: Ania.Otwinowska@math.u-psud.fr
\end{document}